\documentclass[a4paper]{article}
\usepackage{amsfonts,amssymb,amsmath,graphicx,placeins}
\usepackage[a4paper,left=25mm,right=25mm,top=25mm,bottom=30mm]{geometry}
\usepackage[auth-lg,affil-it]{authblk}
\usepackage[unicode,colorlinks=true,urlcolor=black,linkcolor=black,citecolor=black,pdfborder={0 0 0}]{hyperref}
\usepackage{cite}

\begin{document}

\title{Subdivision Shell Elements with Anisotropic Growth}

\author[1]{Roman Vetter}
\author[1,2]{Norbert Stoop}
\author[1]{Thomas Jenni}
\author[1]{Falk K.~Wittel}
\author[1]{Hans J.~Herrmann}
\affil[1]{Computational Physics for Engineering Materials, IfB, ETH Zurich, Schafmattstrasse 6, CH-8093 Zurich, Switzerland}
\affil[2]{IAS Institute of Applied Simulations, ZHAW Zurich University of Applied Sciences, CH-8820 W\"adenswil, Switzerland}

\maketitle

\begin{abstract}
A thin shell finite element approach based on Loop's subdivision surfaces is proposed, capable of dealing with large deformations and anisotropic growth. To this end, the Kirchhoff-Love theory of thin shells is derived and extended to allow for arbitrary in-plane growth. The simplicity and computational efficiency of the subdivision thin shell elements is outstanding, which is demonstrated on a few standard loading benchmarks. With this powerful tool at hand, we demonstrate the broad range of possible applications by numerical solution of several growth scenarios, ranging from the uniform growth of a sphere, to boundary instabilities induced by large anisotropic growth. Finally, it is shown that the problem of a slowly and uniformly growing sheet confined in a fixed hollow sphere is equivalent to the inverse process where a sheet of fixed size is slowly crumpled in a shrinking hollow sphere in the frictionless, quasi-static, elastic limit.
\end{abstract}

\section{Introduction}

Thin shells are vital parts of innumerable problems in structural engineering and material science. It is impossible to imagine many of today's technological achievements without a deep understanding of their large deformation response to external loads, spatial constraints, nonlinear material behavior, self-contact and other multifarious interactions. For solving the complex interplay of these effects, one resorts to numerical methods, like for the wrinkling of metal sheets in a vehicle crash (e.g., \cite{O98,WWK01,MMIKH01,MO05}), or the crumpling of paper \cite{LGLMW95,BAP97,BK05,VG06,TAT08,TAT09,VG11}, just to name a few. The finite element method (FEM) has proven to be among the most flexible and efficient tools for a large number of such problems, in particular where complicated geometries, strong material nonlinearities, or anisotropy come into play.

The nontrivial nature of the various large strain modes of shells is due to their thinness, enabling a large variety of complex three-dimensional deformations induced by the aforementioned external or intrinsic constraints. In the Kirchhoff-Love theory \cite{L88}, which has been well understood and widely applied for decades, the thinness of the shell manifests itself in the assumption that material lines that are straight and perpendicular to the middle surface retain these properties and their length during shell deformations. To account for out-of-plane bending stiffness, the resulting total elastic energy formulation integrates the mean and Gaussian curvatures over the shell's middle surface. In the context of a finite element treatment of structural shell analysis, this inevitably calls for shape functions with continuous first derivatives ($C^1$ continuity) across element boundaries, i.e., functions that belong to the Sobolev space $H^2$. This requirement has posed a tough challenge in the history of shell finite elements. Many of the developed elements introduce interpolation coefficients for higher derivatives of the displacement field, leading to a significant increase in the number of degrees of freedom (DOFs) to compute. Parisch \cite{P95} has proposed quadrilateral shell elements with only displacement variables at the surface nodes and an auxiliary degree of freedom on the middle surface. An entire class of quadrilateral $C^1$ elements is obtained from tensor products of Hermite polynomials, see e.g. Ref.~\cite{ZT05}. However, while simple to construct, they are limited to regular rectangular meshes. The family of Hsieh-Clough-Tocher triangles have shown success in plate bending and other biharmonic problems \cite{SC07}, but are tedious to set up and add many additional DOFs. The recently developed subdivision surface shape functions \cite{COS00,CO01}, which strictly belong to the class $C^1$ everywhere in the domain, have successfully vanquished all these problems, and are hence gaining increasing interest in the FE community.

External loads and constraints are not the only causes for nonlinear shell deformations. Thinking of nature's soft tissues such as leaves, flower petals, cell membranes, insect wings, etc., it becomes evident that large deformations of shells are often on account of growth or shrinkage (e.g., \cite{T95,DBA08}). We use the term \textit{growth} to represent both these mutually inverse processes in the following. Growth often leads to the inevitable development of residual stresses (e.g., \cite{H86,SZJNH96}), causing a shell to deform in an attempt to minimize them. The study of growing thin shells is, however, not limited to bioengineering. There is also a large potential in the blossoming fields of bionics or material engineering, where a variety of smart or self-actuating materials is designed (e.g., \cite{KHBSH12,KHHS12}). Predicting the large deformation response of such growing thin sheets, be it for whatever physical cause or technical purpose, calls for an efficient and robust numerical tool inherently featuring the capability to grow according to arbitrary prescribed anisotropic growth fields. Not many common numerical discretization techniques are fit for this kind of anisotropy. Discrete elements, or beam networks, for instance, are not well suited, as their only degree of freedom capable of accounting for in-plane growth is the edge length connecting the mesh vertices, which may not be aligned with the desired growth direction. To the finite element method, on the other hand, material anisotropy poses no problem, since the mesh faces are numerically integrated over, regardless of their orientation.

The purpose of this article is to present easy-to-implement thin shell finite elements that are capable of anisotropic in-plane growth, while profiting from the superior efficiency and strong robustness of the subdivision surface paradigm at the same time. The adopted continuum model is based on the geometrically nonlinear Kirchhoff-Love shell theory, which we extend by a large-strain continuum growth model that supports anisotropic growth fields. The basic idea is to assume that body deformations can be due to both a change of mass or volume and an elastic response \cite{H68,SZJNH96}. The thin shell growth model employed here is based on the volumetric growth assumptions of Rodriguez et al. \cite{RHM94}, which have been put on rigorous foundation \cite{DQ02,LH02}.

The paper is organized as follows: The next section extends the Kirchhoff-Love theory of thin shells for large deformations by Rodriguez' growth ansatz, for the continuum shell and in discretized form suitable to finite element analysis. A short overview of subdivision surface interpolation for finite elements and boundary conditions is given in the subsequent section \ref{sec:subdiv}, followed by a description of some implementation details in section \ref{sec:impl}. Standard load cases to verify and benchmark our shell elements with and without growth are presented in the last section, where we also demonstrate the huge potential in applicability of growing thin shells to several problems in material science and engineering.

\section{The Kirchhoff-Love Theory with in-Plane Growth}

In the following section, the Kirchhoff-Love theory is briefly derived following the common \textit{stress-resultant} formulation, meaning that the stresses are integrated analytically through the thickness, so that one is left with a resultant stress on the middle-surface. In the course, the theory is amended by anisotropic in-plane growth for large strains.

\subsection{Kinematics of Deformation}

Let Greek indices $\alpha,\beta,\gamma,\delta$ take the values $1$ and $2$, and Latin indices $i,j$ take values from $1$ to $3$. Unless otherwise mentioned, the Einstein summation convention applies to repeated indices, and lower (upper) indices denote the covariant (contravariant) components.

\begin{figure*}[htpb]
	\begin{center}
	\includegraphics{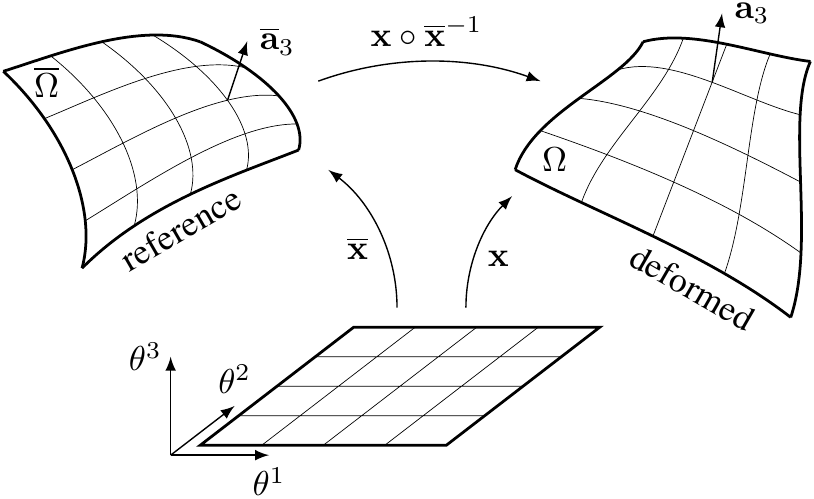}
	\caption{Reference and deformed configurations of the shell middle surface with parameterizations ${\overline{\bf x}}(\theta^1,\theta^2)$ and ${\bf x}(\theta^1,\theta^2)$, respectively.}
	\label{fig:coord_sys}
	\end{center}
\end{figure*}

Let $\overline{\Omega} \subset \mathbb{E}^3$ be the geometry of the stress-free undeformed middle surface of a shell with small thickness $h$, embedded in three-dimensional Euclidean space. Under the action of external loads or growth, the shell deforms into a new configuration characterized by the middle surface $\Omega \subset \mathbb{E}^3$. Let $\{\theta^1,\theta^2, \theta^3\}$ be a curvilinear coordinate system, and $\overline{{\bf x}}(\theta^1,\theta^2)$ and ${\bf x}(\theta^1,\theta^2)$ be parametrizations of the reference middle surface $\overline{\Omega}$ and deformed middle surface $\Omega$, respectively, see Fig.~\ref{fig:coord_sys}. The positions $\overline{{\bf r}}$ and ${\bf r}$ of material points in the reference and deformed shell may be parametrized as
\begin{align}
\overline{{\bf r}}(\theta^1,\theta^2,\theta^3) &= \overline{{\bf x}}(\theta^1,\theta^2) + \theta^3 \overline{{\bf a}}_3(\theta^1,\theta^2),\qquad\theta^3 \in [-h/2,h/2],\\
\label{eq:r_def}
{\bf r}(\theta^1,\theta^2,\theta^3) &= {\bf x}(\theta^1,\theta^2) + \theta^3 {\bf a}_3(\theta^1,\theta^2),\qquad\theta^3 \in [-h/2,h/2].
\end{align}
The tangent space of the middle surface is spanned by the vectors
\begin{equation}
\overline{{\bf a}}_{\alpha}(\theta^1,\theta^2) = \frac{\partial\overline{\bf x}}{\partial\theta^{\alpha}} =: \overline{{\bf x}}_{,\alpha} \qquad \textrm{and} \qquad {\bf a}_{\alpha}(\theta^1,\theta^2) = \frac{\partial\bf x}{\partial\theta^{\alpha}} =: {\bf x}_{,\alpha},
\end{equation}
and by virtue of the Kirchhoff kinematic assumption, the material orientation in the thickness direction of the shell is determined by the \textit{shell directors}
\begin{equation}
\label{eq:conda3}
{\bf a}_3 = \frac{{\bf a}_1\times {\bf a}_2}{|{\bf a}_1\times {\bf a}_2|}\qquad\textrm{and}\qquad{\overline{\bf a}}_3 = \frac{{\overline{\bf a}}_1\times {\overline{\bf a}}_2}{|{\overline{\bf a}}_1\times {\overline{\bf a}}_2|}.
\end{equation}

The covariant components of the surface metric tensors, or first fundamental forms, follow as
\begin{equation}
\overline{a}_{\alpha\beta} = \overline{{\bf a}}_{\alpha} \cdot \overline{{\bf a}}_{\beta},
\qquad
a_{\alpha\beta} = {\bf a}_{\alpha} \cdot {\bf a}_{\beta},
\end{equation}
those of the shape tensors, or second fundamental forms, as
\begin{equation}
\overline{\kappa}_{\alpha\beta} = - \overline{{\bf a}}_{3,\alpha} \cdot \overline{{\bf a}}_{\beta} = \overline{{\bf a}}_{3} \cdot \overline{{\bf a}}_{\alpha,\beta},
\qquad
\kappa_{\alpha\beta} = - {\bf a}_{3,\alpha} \cdot {\bf a}_{\beta} = {\bf a}_{3} \cdot {\bf a}_{\alpha,\beta},
\end{equation}
and the infinitesimal area element can be expressed as $\textrm{d}\overline{\Omega} = |\overline{{\bf a}}_1 \times \overline{{\bf a}}_2|\;\textrm{d}\theta^1 \textrm{d}\theta^2$. The covariant basis vectors for a generic point within the shell are given by
\begin{align}
\overline{{\bf g}}_{\alpha} = \frac{\partial\overline{{\bf r}}}{\partial\theta^{\alpha}} = \overline{{\bf a}}_{\alpha} + \theta^3 \overline{{\bf a}}_{3,\alpha}&,\qquad\overline{{\bf g}}_{3} = \frac{\partial\overline{{\bf r}}}{\partial\theta^{3}} = \overline{{\bf a}}_{3},\\ 
{\bf g}_{\alpha} = \frac{\partial{\bf r}}{\partial\theta^{\alpha}} = {\bf a}_{\alpha} + \theta^3 {\bf a}_{3,\alpha}&,\qquad{\bf g}_{3} = \frac{\partial{\bf r}}{\partial\theta^{3}} = {\bf a}_{3},
\end{align}
and the covariant components of the corresponding metric tensors $\overline{g}$, $g$ are
\begin{equation}
\overline{g}_{ij}= \overline{{\bf g}}_i \cdot \overline{{\bf g}}_j,\qquad g_{ij} = {\bf g}_i \cdot {\bf g}_j.
\end{equation}
Owing to the Kirchhoff constraints, Eq.~(\ref{eq:conda3}), we require that $g_{33}=1$ and $g_{\alpha 3} = g_{3 \alpha}= 0$. The deformation gradient in curvilinear coordinates reads
\begin{equation}
{\bf F} := \nabla_{\overline{\bf r}}{\bf r} = \frac{\partial {\bf r}}{\partial \theta^i} \otimes \overline{{\bf g}}^i =  {\bf g}_{i} \otimes \overline{{\bf g}}^i
\end{equation}
and maps between deformed and reference metric,
\begin{equation}\label{eq:fmap}
{\bf F}\overline{{\bf g}}_j = {\bf g}_i \delta^i_j = {\bf g}_j.
\end{equation}
The Green-Lagrange strain tensor ${\bf E}$ in curvilinear coordinates is then defined as
\begin{equation}
{\bf E} := \frac{1}{2}(g-\overline{g}) = \frac{1}{2} ({\bf F}^{\textrm{T}} \overline{g} {\bf F} - \overline{g}).
\end{equation}
In the growth model proposed by Rodriguez et al. \cite{RHM94}, the geometric deformation gradient is multiplicatively decomposed into a growth tensor ${\bf G}$ and a purely elastic response ${\bf A}$, that ensures compatibility and continuity of the body, according to
\begin{equation}
{\bf F} = {\bf A\,G},
\end{equation}
analogous to Lee's multiplicative decomposition in elastic-plastic modeling \cite{L69}. If we require that ${\bf G}$ be independent of the deformed configuration, the strains due to the elastic response ${\bf A}$ become
\begin{equation}
\label{eq:el_strain}
{\bf E}_{\textrm{e}} = \frac{1}{2} ({\bf A}^{\textrm{T}} \overline{g} {\bf A} - \overline{g}) = \frac{1}{2} ({\bf G}^{-\textrm{T}} {\bf F}^{\textrm{T}} \overline{g} {\bf F} {\bf G}^{-1} - \overline{g}) = \frac{1}{2} ({\bf G}^{-\textrm{T}} g {\bf G}^{-1} - \overline{g}) = \frac{1}{2}(\tilde{g}-\overline{g}),
\end{equation}
where $\tilde{g} := {\bf G}^{-\textrm{T}}g{\bf G}^{-1}$ is the growth-modified metric. For the following, we restrict ${\bf G}$ to anisotropic in-plane growth, i.e.,
\begin{equation}
{\bf G} = [G_{ij}] = \begin{bmatrix}[G_{\alpha\beta}] & 0\\0^{\textrm{T}} & 1\end{bmatrix}.
\end{equation}
Consequently, $\tilde{g}_{33}=1$ and $\tilde{g}_{\alpha 3} = \tilde{g}_{3 \alpha}= 0$, and the remaining components can be expanded in terms of the thickness parameter $\theta^3$ to
\begin{align}
\tilde{g}_{\alpha \beta} &= \left([G_{\alpha \beta}]^{-\textrm{T}}\,[a_{\alpha \beta}]\,[G_{\alpha \beta}]^{-1}\right)_{\alpha \beta} - 2 \theta^3 \left([G_{\alpha \beta}]^{-\textrm{T}}\,[\kappa_{\alpha \beta}]\,[G_{\alpha \beta}]^{-1}\right)_{\alpha \beta} + \mathcal{O}((\theta^3)^2)\\
&= \tilde{a}_{\alpha \beta} - 2 \theta^3 \tilde{\kappa}_{\alpha \beta} + \mathcal{O}((\theta^3)^2),
\end{align}
Analogous to the derivation of the Kirchhoff shell without growth, we neglect terms $\mathcal{O}((\theta^3)^2)$ in the following, making the theory a \textit{first order shell theory} which is valid for small thicknesses $h$. The non-zero components of the elastic strain tensor (\ref{eq:el_strain}) then follow as
\begin{equation}
E_{\alpha\beta} \approx \tilde{\alpha}_{\alpha\beta} + \theta^3 \tilde{\beta}_{\alpha\beta}
\end{equation}
where $\tilde{\alpha}$ and $\tilde{\beta}$ are the growth-modified membrane and bending strain tensors, respectively,
\begin{align}
\tilde{\alpha}_{\alpha \beta} &= \frac{1}{2}(\tilde{a}_{\alpha \beta} - \overline{a}_{\alpha \beta})\\
\tilde{\beta}_{\alpha \beta} &= \overline{\kappa}_{\alpha \beta} - \tilde{\kappa}_{\alpha \beta}.
\end{align}

The above expressions are formally identical to the strains of a Kirchhoff shell without growth. Note also that the four in-plane growth parameters $G_{\alpha\beta}$, which in general may be functions of various external or internal variables such as time, space, stress, etc. \cite{AG05,MG11}, are expressed with respect to the reference tangent basis $\{\overline{\bf a}_1,\overline{\bf a}_2\}$ in the above formalism. In practice, it is thus necessary to perform a change of basis when they are to be given with respect to a specific coordinate system, such as the Cartesian, for instance.

\subsection{Constitutive Model}

Assuming that the shell obeys the St.~Venant-Kirchhoff law, the connection between its geometrical configuration and material properties is provided by the Koiter energy density functional \cite{K66,C00}
\begin{equation}
\label{eq:koiter}
W = \frac{1}{2} K H^{\alpha\beta\gamma\delta}\tilde{\alpha}_{\alpha\beta} \tilde{\alpha}_{\gamma\delta} + \frac{1}{2} D H^{\alpha\beta\gamma\delta}\tilde{\beta}_{\alpha\beta} \tilde{\beta}_{\gamma\delta}
\end{equation}
with membrane stiffness $K$ and bending rigidity $D$, given by
\begin{equation}
K=\frac{Yh}{1-\nu^2},\qquad D=\frac{Yh^3}{12(1-\nu^2)},
\end{equation}
and with the elasticity tensor
\begin{align}
H^{\alpha\beta\gamma\delta} = \nu\overline{a}^{\alpha\beta}\overline{a}^{\gamma\delta}+\frac{1-\nu}{2}(\overline{a}^{\alpha\gamma}\overline{a}^{\beta\delta}+\overline{a}^{\alpha\delta}\overline{a}^{\beta\gamma}).
\end{align}
$Y$ is the Young's modulus and $\nu$ the Poisson ratio. The growth-modified resultant membrane stresses $\tilde{n}$ and resultant bending stresses $\tilde{m}$ follow by the principle of work conjugacy:
\begin{equation}
\tilde{n}^{\alpha\beta} = \frac{\partial W}{\partial\tilde{\alpha}_{\alpha\beta}} = K H^{\alpha\beta\gamma\delta}\tilde{\alpha}_{\gamma\delta},\qquad \tilde{m}^{\alpha\beta} = \frac{\partial W}{\partial\tilde{\beta}_{\alpha\beta}} = D H^{\alpha\beta\gamma\delta}\tilde{\beta}_{\gamma\delta}. 
\end{equation}

Notice that the energy density in the above form coincides with the non-Euclidean plate approach derived in Ref.~\cite{ESK09}.

\subsection{Variational Formulation}

The total potential energy $\Phi$ of the Koiter shell with total Lagrangian displacement of the middle surface, ${\bf u} = {\bf x} - \overline{\bf x}$, is obtained by adding the internal elastic energy $\Phi^{\textrm{int}}$ to the contribution $\Phi^{\textrm{ext}}$ from external loads ${\bf q}$ per unit surface area and traction ${\bf N}$ per unit edge length, yielding
\begin{align}
\label{eq:totalenergy}
\Phi[{\bf u}] &= \Phi^{\textrm{int}}[{\bf u}] + \Phi^{\textrm{ext}}[{\bf u}],\\
\Phi^{\textrm{int}}[{\bf u}] &= \int_{\overline{\Omega}} W\,\textrm{d}\overline{\Omega},\\
\Phi^{\textrm{ext}}[{\bf u}] &= -\int_{\overline{\Omega}}{\bf q} \cdot {\bf u} \,\textrm{d}\overline{\Omega}  - \int_{\partial{\overline{\Omega}}}{\bf N} \cdot {\bf u}\,\textrm{d}\overline{s}.
\end{align}

Our aim is to find the minimum of Eq. (\ref{eq:totalenergy}) for prescribed growth tensors, which is equivalent to finding a displacement field ${\bf u}$ satisfying the variational problem
\begin{align}
\label{eq:weak_formulation}
0 &= \delta \Phi[{\bf u}] = \delta \Phi^{\textrm{int}}[{\bf u}] + \delta \Phi^{\textrm{ext}}[{\bf u}],\\
\delta \Phi^{\textrm{int}}[{\bf u}] &= \int_{\overline{\Omega}} \left( \tilde{n}^{\alpha\beta}\,\delta\tilde{\alpha}_{\alpha\beta} + \tilde{m}^{\alpha\beta}\,\delta\tilde{\beta}_{\alpha\beta} \right)\,\textrm{d}\overline{\Omega},\\
\delta \Phi^{\textrm{ext}}[{\bf u}] &= -\int_{\overline{\Omega}} {\bf q}\cdot \delta {\bf u}\,\textrm{d}\overline{\Omega} - \int_{\partial\overline{\Omega}} {\bf N} \cdot \delta {\bf u}\,\textrm{d}\overline{s}.
\end{align}
We further augment Eq.~(\ref{eq:weak_formulation}) with the usual inertial term to capture the dynamics of the system. The variational statement thus becomes
\begin{equation}
\label{eq:dyn_weak_formulation}
0 = \delta \Phi[{\bf u}] + \int_{\overline{\Omega}} h \rho \ddot{{\bf u}} \cdot \delta {\bf u}\,\textrm{d}\overline{\Omega},
\end{equation}
where $\rho$ is the mass density of the shell. The variation of the membrane and bending strains are easily calculated if the growth tensor is assumed not to depend on the displacement field:
\begin{equation}
\label{eq:strain_variation}
[\delta \tilde{\alpha}_{\alpha \beta}] = [G_{\alpha\beta}]^{-\textrm{T}}\,[\delta a_{\alpha \beta}]\,[G_{\alpha\beta}]^{-1},
\qquad
[\delta \tilde{\beta}_{\alpha \beta}] = [G_{\alpha\beta}]^{-\textrm{T}}\,[\delta \kappa_{\alpha \beta}]\,[G_{\alpha\beta}]^{-1}.
\end{equation}
$\delta a_{\alpha \beta}$ and $\delta \kappa_{\alpha \beta}$ are the usual first strain variations for shells without growth, see e.g. Ref.~\cite{SF89,SFR90}. Since we will integrate the Newtonian equations of motion induced by Eq. (\ref{eq:dyn_weak_formulation}) explicitly in time (cf.~Section \ref{sec:time}), no second variations are needed at this point. We note, however, that the derivation of the second variations straightforwardly follows the usual formalism (e.g., \cite{SFR90,G03}) without complications for our growth-modified strains.

\subsection{Finite Element Discretization}

The finite element discretization is done in the usual way: The minimization problem is replaced by an approximate minimization problem over a finite subspace $V_h \subset V := H^2(\overline{\Omega},\mathbb{E}^3)$ of admissible displacements:
\begin{equation}
\label{minimization_problem}
\inf\limits_{{\bf u}\in V}\Phi[{\bf u}]\quad\longrightarrow\quad\min\limits_{{\bf u}_h\in V_h} \Phi[{\bf u}_h].
\end{equation}
$V_h$ is spanned by a finite set of basis functions $\{N_I(\theta^1,\theta^2), I = 1,..., \mathcal{N}_{\textrm{n}}\}$ with local support, where $\mathcal{N}_\textrm{n}$ is the number of mesh nodes. The displacement field is then written as a linear combination of the trial space basis functions:
\begin{equation}
\label{shell_displacement_field_interpolation}
{\bf u}_h(\theta^1,\theta^2) = \sum_{I=1}^{\mathcal{N}_\textrm{n}}{\bf u}_I N_I(\theta^1,\theta^2),
\qquad
\delta{\bf u}_h(\theta^1,\theta^2) = \sum_{I=1}^{\mathcal{N}_\textrm{n}}\delta{\bf u}_I N_I(\theta^1,\theta^2).
\end{equation}
Substituting the above interpolation into the weak form, Eq.~(\ref{eq:dyn_weak_formulation}), and using the arbitrariness of the trial field, the variational statement is recast into an algebraic minimization problem for the nodal displacements ${\bf u}_I$:
\begin{align}
\label{shell_generalized_forces}
0 &= {\bf f}^{\textrm{int}}_I - {\bf f}^{\textrm{ext}}_I + \sum_J M_{IJ}\ddot{{\bf u}}_J,\\
\label{eq:fint}
{\bf f}^{\textrm{int}}_I &= - \int_{\overline{\Omega}} \left(\tilde{n}^{\alpha\beta}\frac{\partial \tilde{\alpha}_{\alpha\beta}}{\partial {\bf u}_I} + \tilde{m}^{\alpha\beta}\frac{\partial \tilde{\beta}_{\alpha\beta}}{\partial {\bf u}_I} \right)\,\textrm{d}\overline{\Omega},\\
\label{eq:fext}
{\bf f}^{\textrm{ext}}_I &= \int_{\overline{\Omega}} {\bf q}N_I\,\textrm{d}\overline{\Omega} + \int_{\partial\overline{\Omega}} {\bf N} N_I\,\textrm{d}\overline{s},\\
\label{eq:mass_integral}
M_{IJ} &= \int_{\overline{\Omega}} h \rho N_I N_J\,\textrm{d}\overline{\Omega}.
\end{align}

As usual in finite element analysis, the integrals (\ref{eq:fint}--\ref{eq:mass_integral}) are evaluated numerically using a quadrature rule with $\mathcal{N}_{\textrm{p}}$ integration points $\{q_p = (\theta^1_p,\theta^2_p), p = 1,..., \mathcal{N}_{\textrm{p}}\}$ and corresponding weights $\{w_p, p = 1,..., \mathcal{N}_{\textrm{p}}\}$, taking advantage of the local support of the shape functions. A single element's contribution to the generalized internal force (\ref{eq:fint}), for instance, becomes
\begin{equation}
\label{eq:fint_elem}
{\bf f}^{\textrm{int}}_{I,e} = - \sum\limits_{p=1}^{\mathcal{N}_{\textrm{p}}} \left[\left(\tilde{n}^{\alpha\beta}\frac{\partial \tilde{\alpha}_{\alpha\beta}}{\partial {\bf u}_I} + \tilde{m}^{\alpha\beta}\frac{\partial \tilde{\beta}_{\alpha\beta}}{\partial {\bf u}_I}\right)|\overline{{\bf a}}_1 \times \overline{{\bf a}}_2|\right]_{(\theta^1_p,\theta^2_p)}w_p,
\end{equation}
where $[\,\cdot\,]_{(\theta^1_p,\theta^2_p)}$ denotes evaluation of the integrand at the quadrature point $q_p$ mapped onto element $e$. Analogous to Eq.~(\ref{eq:strain_variation}), the growth-modified strain derivatives are straightforward to calculate since the growth tensor is independent of the deformed configuration:
\begin{equation}
\left[\frac{\partial \tilde{\alpha}_{\alpha\beta}}{\partial {\bf u}_I}\right] = [G_{\alpha\beta}]^{-\textrm{T}}\left[\frac{\partial\alpha_{\alpha\beta}}{\partial {\bf u}_I}\right][G_{\alpha\beta}]^{-1},
\qquad
\left[\frac{\partial \tilde{\beta}_{\alpha\beta}}{\partial {\bf u}_I}\right] = [G_{\alpha\beta}]^{-\textrm{T}}\left[\frac{\partial\beta_{\alpha\beta}}{\partial {\bf u}_I}\right][G_{\alpha\beta}]^{-1}.
\end{equation}

The energy functional (\ref{eq:koiter}) and the generalized internal force (\ref{eq:fint}) contain second derivatives of the displacement field ${\bf u}$. For boundedness of these integrals, the shape functions $N_I$ therefore need $C^1$ continuity. Loop subdivision surfaces, being $C^1$-continuous everywhere and even $C^2$-continuous except at a finite set of extraordinary points, fully meet this requirement.

A noteworthy difference of the above finite element description of growth to recent approaches via the introduction of prescribed non-Euclidean target metrics \cite{MSSR03,AB03,KES07,MDS07} is that the numerical implementation of the present model includes the change of reference curvature when the surface grows, which is missing in the tethered mass-spring model of the target metric approach \cite{MSSR03,MDS07}.

\section{Subdivision Surfaces}
\label{sec:subdiv}

Subdivision surfaces were developed simultaneously by Catmull and Clark \cite{CC78} in the context of computer graphics in 1978 as a method of representing smooth surfaces by a coarse polygonal mesh, termed the \textit{control mesh}. Cirak et al.~\cite{COS00,CO01} have ported them to the finite element method.

\subsection{Shape Functions}

The methodology is based on Stam's eigenanalysis \cite{S99} of Loop's recursive refinement rule \cite{L87} for triangulated surfaces with arbitrary topology, which gave access to a set of 12 quartic box splines, that exactly interpolate the infinitely refined surface, called \textit{limit surface}, at all points except a finite set. In fundamental difference to traditional finite elements, subdivision surfaces gain $C^1$ continuity, or $H^2$ integrability, required for a finite deformation energy of a shell, at the expense of a larger local support of the basis functions. Instead of only the \textit{1-ring} consisting of directly adjacent neighbor elements, each element's support spans also the \textit{2-ring} of next-neighboring elements. Details on the Loop subdivision shape functions and their application to interpolating the limit surface on arbitrarily triangulated meshes can be found in Refs.~\cite{S99,COS00}.

Subdivision surface elements are much more efficient in explicitly integrated Newton's dynamics than conventional $C^0$ shell elements, as will be demonstrated toward the end of this paper. They go without rotational degrees of freedom, i.e., each mesh node carries three displacement variables only, which reduces the overall system size by a factor of two or more when compared to traditional three-node thin shell elements. Moreover, a single Gauss point per element has been found sufficient for convergence and accuracy in previous linear \cite{COS00} and nonlinear models \cite{CO01}. This apparently stems from the enhanced support of the shape functions. While classical triangular elements with quartic polynomials require at least six integration points per triangle \cite{S71}, each element is integrated using only its own points. The numerical integration over subdivision surface elements, on the other hand, includes all points from the element's 1-ring (typically 12), so a single point per triangle suffices to satisfy the theoretical lower bound. Indeed, we have not found any significant inaccuracies or spurious modes when using a one-point quadrature even in situations with extremely large deformation (see Section \ref{sec:ex}). Of course, increasing the number of quadrature points may assist in resolving material or growth anisotropies.

On the downside, the extended support of subdivision shape functions does require new concepts for boundary conditions, as explained in the next subsection. Also, the stiffness bandwidth is increased by a factor of up to two on average compared to traditional plain triangular elements with six DOFs per mesh node, compromising somewhat the gain from avoiding rotational DOFs in the first place in analyses where the stiffness matrix is actually assembled. In transient analysis with explicit integration in time, however, where no linear system of equations needs to be solved, the overall computational advantage clearly prevails.

\subsection{Boundary Conditions}

The treatment of domain boundaries is in general non-trivial for subdivision surface interpolation because elements along a shell boundary lack a complete 1-ring and hence \textit{per se} cannot be interpolated like elements in the interior of the domain. Three principal approaches have been proposed to solve this problem.

\begin{figure*}[htpb]
	\begin{center}
	\includegraphics{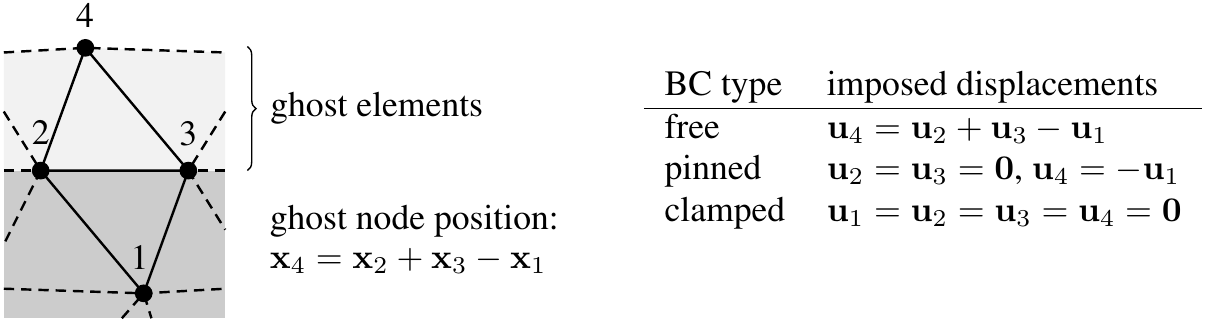}
	\caption{Boundary conditions of the Schweitzer-Cirak type.}
	\label{fig:bc}
	\end{center}
\end{figure*}

\begin{enumerate}
\setlength{\itemsep}{0pt}
\setlength{\parskip}{0pt}
\item As pointed out by Cirak et al. and Biermann et al.~\cite{COS00,CL11,BLZ00}, a special subdivision rule may be applied to boundary elements, corresponding to one-dimensional subdivision on the boundary curve.

\item Alternatively, Schweitzer \cite{S96} suggested appending an additional row of ``ghost'' vertices and elements to the control mesh along its boundaries. If these ghost vertices are positioned and the displacement field is projected onto them component-wise according to Cirak et al.~\cite{COS00}, see Fig.~\ref{fig:bc}, the usual boundary conditions can easily be imposed. We will refer to this type of boundary as the Schweitzer-Cirak type in the following. However, as demonstrated by Green \cite{G03,G04}, such boundary constraints are overly restrictive and lead to a reduction of the convergence with the number of elements from order two to one. In practice, this implies that rather fine meshes are required for high accuracy, undermining the otherwise high computational efficiency of subdivision finite elements.

\item The third approach was proposed by Green \cite{G03,G04} as a remedy to the limitations of the second. Instead of drastically constraining the ghost displacements according to Fig.~\ref{fig:bc}, only the minimum set of necessary conditions is imposed directly on the limit surface. The resulting linear constraint equations can be solved using the penalty method, Lagrange multipliers, or any other solving technique suitable for constrained minimization.
\end{enumerate}

In summary, Schweitzer-Cirak boundaries are the easiest to implement, but should be avoided in cases where physical accuracy is crucial. The two alternatives are significantly more complex to implement in general, with the exception of free Green boundaries, where the ghost nodes are simply left unconstrained. Such ghost nodes are, however, unsuited for static analysis due to underdetermination of the system. Most boundaries in the numerical examples presented in Section \ref{sec:ex} are free, and since we solve the dynamic problem, Green's convenient method is applied there. We use Schweitzer-Cirak boundaries in all other situations for simplicity.

\section{Numerical Implementation}
\label{sec:impl}

Thin shells are commonly formulated and implemented using Voigt's vector notation for convenience and efficiency, exploiting the symmetry of the involved tensors. Readers interested in the details are referred to the numerous publications describing the technicalities, such as Refs.~\cite{SFR89,SFR90,COS00}. The same principles apply straightforwardly to the growth-modified strains and stress resultants. In this section, we focus on the integration of dynamics and shell contacts.

\subsection{Integration in Time}
\label{sec:time}

We solve the hyperbolic equilibrium equations
\begin{equation}
{\bf M}\ddot{{\bf u}}(t) + {\bf C}{\dot{{\bf u}}}(t) + {\bf f}^{\textrm{int}}\left({\bf u}(t)\right) = {\bf f}^{\textrm{ext}}\left({\bf u}(t),t\right).
\end{equation}
Here, ${\bf u}(t)$ is the displacement vector containing the nodal interpolants ${\bf u}_I$ at time $t$, ${\bf M}$ is the mass matrix assembled from element contributions to Eq.~(\ref{eq:mass_integral}), and ${\bf C}$ is a viscous damping matrix for equilibration. ${\bf f}^{\textrm{int}}$ and ${\bf f}^{\textrm{ext}}$ account for internal out-of-equilibrium forces, Eq.~(\ref{eq:fint_elem}), and external loads and contact forces, respectively. Newmark's family of integration methods \cite{N59} is widely used in structural dynamics \cite{ZWLX92,R04}. Let ${\bf u}_t \approx {\bf u}(t)$, ${\bf v}_t \approx \dot{\bf u}(t)$ and ${\bf a}_t \approx \ddot{\bf u}(t)$ be the discretized approximations of the displacement vector and its time derivatives. For fixed Newmark scheme parameters $\beta$ and $\gamma$, they are integrated according to
\begin{align}
{\bf u}_{t+\Delta t} &= {\bf u}_t + \Delta t\,{\bf v}_t + \frac{(\Delta t)^2}{2}\big((1-2\beta){\bf a}_t + 2\beta\,{\bf a}_{t+\Delta t}\big),\\
{\bf v}_{t+\Delta t} &= {\bf v}_t + \Delta t\big((1-\gamma){\bf a}_t + \gamma\,{\bf a}_{t+\Delta t}\big).
\end{align}
$\Delta t$ denotes the finite time step. We apply the unconditionally stable \textit{constant-average acceleration method}, that is obtained by setting $\beta=1/4$ and $\gamma=1/2$, in form of a predictor-corrector scheme with lumped masses and subcritical lumped damping. For adaptive time step control, the \textit{a posteriori} local error estimator by Zienkiewicz and Xie \cite{ZX91,ZWLX92} is employed.

\subsection{Shell-Shell Contact}

For problems where the shell may be in contact with itself, we use a hierarchical spatial decomposition to find contact points efficiently, similarly to Ref.~\cite{GCSO99}. On the coarsest level, elements are placed into an array of cubic axis-aligned cells. Each cell holds a reference to all elements whose axis-aligned boundary boxes (AABBs), extended by half the shell thickness in each direction, overlap with it. Each element pair sharing at least one cell enters the mid-level check where the elements' extended AABBs are tested for overlap. If they do, the distance between the two elements and the corresponding closes points are computed on the finest level of the hierarchy.

Finding the closest points between the subdivision limit surfaces of two triangles is a four-dimensional nonlinear optimization problem with linear inequality constraints. As such, it may be solved with, e.g., the projected gradient method \cite{R60}, a Newton-type penalty method, or any other scheme dedicated to such problems. The limit surface can be evaluated at any point using the ideas of Stam \cite{S99}. To reduce the computational expenses, we apply collision detection on the faceted control mesh instead. There, the problem boils down to evaluating the distance between all nine edge-edge pairs and six vertex-face pairs belonging to the two triangular elements in question \cite{BFA02,E05}. As a compromise, one could consider doing this on any level of subdivision refinement, rather than on the control mesh, like has been done in Ref.~\cite{GCSO99}.

The closest points on edge-edge pairs are efficiently determined using an algorithm by Sunday \cite{S01} with edges normalized to unit length. Eriscon's robust algorithm \cite{E05} is used for the vertex-triangle pairs. Once the closest points have been identified on both contacting triangles, and if they are less than the shell thickness apart, a repulsive contact force that is proportional to the surface area of the contacting elements is distributed to the six involved vertices using linear interpolation. The precise form of the contact law turned out to be irrelevant in many practical applications. An important feature is divergence at zero distance to avoid interpenetration when the shell is very thin.

\section{Numerical Examples}
\label{sec:ex}

A standard verification obstacle course for small shell deformations in the linear regime can be found in Ref.~\cite{COS00}. Although we have verified our implementations using such examples, we are omitting the details here, as our focus is on large deformations with geometric nonlinearity, and growth. A few tests of nonlinear subdivision shell elements without growth can be found in Refs.~\cite{CO01,CL11} for different constitutive models. In the following, we perform a couple of standard tests to verify our implementation of the geometrically nonlinear Koiter shell theory, followed by some examples including various types of growth-induced nonlinearity, with increasing sophistication. The numerical data presented below has been obtained using one quadrature point per element (at the barycenter) unless indicated otherwise. All standard benchmarks were also repeated using a six-point rule to confirm that a single point is indeed sufficient. In all cases, the error from using the one-point rule is of the order of the discretization error or lower. Numerical upper bounds for this relative error are given for each of the standard tests. We denote the subdivision shell finite elements with one-point integration by SD3R in the following for brevity, in the spirit of the Abaqus FEA software \cite{AUM}.

A problem we are not addressing here in great detail is numerical locking. The presented shell model, being an extension of the Koiter shell, obviously avoids \textit{shear locking} with subdivision surface shape functions. \textit{Membrane locking}, however, which is the deterioration of convergence in the zero thickness limit in settings with non-inhibited pure bending, has been shown to exists independently of the smoothness of the approximation space \cite{CB98,CB11}. Indeed, we have found membrane locking with subdivision surface elements in numerical benchmarks designed for detecting locking \cite{CB98,BIC00}. Our numerical tests show that membrane locking becomes prominent in length-to-thickness ratios of $L/h \gtrsim 10^3$. The numerical examples shown below are all either far away from the vanishing thickness limit or inhibit pure bending, such that effects from membrane locking are either very small or completely absent.

\subsection{Inflation and Isotropic Growth of a Sphere}

Only few geometrically nonlinear problems are amenable to analytical solution. The inflation of a sphere is one of them \cite{GS50}. For simplicity, we set the bending rigidity $D=0$ in this example, i.e., a change in energy is assumed purely due to stretching. Since a spherical shell has no boundary, this example is perfect for verifying both the pure response to large membrane stresses, and uniform in-plane growth, in a single scenario. Consider a growth tensor
\begin{equation}
{\bf G} = \textrm{diag}\big(1+g,1+g,1\big)
\end{equation}
with respect to the local tangent basis $\{{\bf a}_1/|{\bf a}_1|,{\bf a}_2/|{\bf a}_2|,{\bf a}_3\}$, where $g$ is a positive growth factor. The sphere with initial radius ${\overline R}$ is then trivially expected to grow uniformly according to $R/{\overline R} =: \lambda = 1+g$. On the other hand, in the absence of bifurcations away from the spherical symmetry \cite{N77}, the pressure $p$ needed to inflate a sphere obeying the Koiter energy density $W$ (\ref{eq:koiter}) from radius ${\overline R}$ to $R\ge{\overline R}$ is easily found by balancing internal and external forces:
\begin{equation}
p = \frac{\partial W}{\partial R}.
\end{equation}
$W=W(R;{\overline R})$ is found using local symmetry on the Green strains
\begin{equation}
E_{\alpha\beta}=\frac{1}{2}(\lambda^2-1)\delta_\alpha^\beta,
\end{equation}
where $\lambda=R/{\overline R}\ge 1$ is the principal stretch. The energy density of the inflated spherical membrane is thus
\begin{equation}
W(R;{\overline R}) = K\frac{(1+\nu)}{4}\left(\lambda^2-1\right)^2,
\end{equation}
yielding the pressure relation
\begin{equation}
p = \frac{Yh}{{\overline R}(1-\nu)}(\lambda^3-\lambda)\ge 0.
\end{equation}

\begin{figure*}[htpb]
	\begin{center}
		\includegraphics{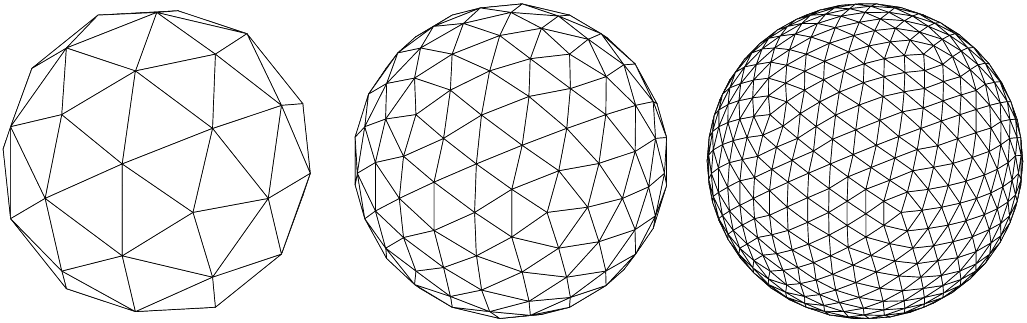}
		\caption{Icosa-spherical meshes for the inflated and growing spherical shell with 80, 320 and 1280 triangles (126, 486 and 1926 DOFs).}
		\label{fig:icosaspheres}
	\end{center}
\end{figure*}

The meshes used for this example are shown in Fig.~\ref{fig:icosaspheres}. They are constructed by recursive quadrisection of the faces of a \textit{regular icosahedron}, followed by a radial projection of the newly created vertices onto the bounding sphere on each level of recursion. The employed recursion depths are 1, 2 and 3, yielding triangulated spheres with 80, 320 and 1280 equilateral elements, respectively. The relevant simulation parameters are ${\overline R}=1$, $h=10^{-3}$, $Y=1$, $\nu=0.3$.

\begin{figure*}[htpb]
	\begin{center}
		\includegraphics{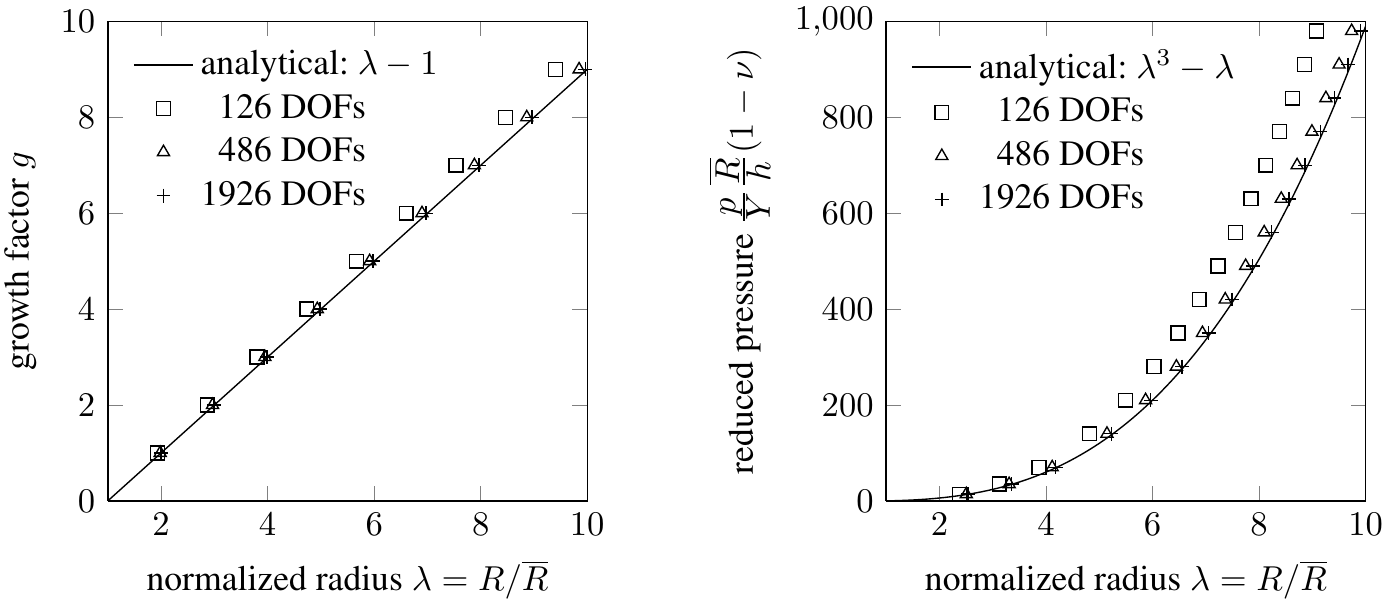}
		\caption{Change of radius of a spherical shell due to isotropic in-plane growth (left) and uniform pressure inflation (right).}
		\label{fig:sphere_pressure}
	\end{center}
\end{figure*}

In Fig.~\ref{fig:sphere_pressure}, we plot the growth-expansion and pressure-expansion curves next to each other, demonstrating the high accuracy and convergence to the analytical solutions in both cases. Relatively small systems are sufficient for high accuracy even at extremely large deformations. The relative error between one-point and six-point quadrature is of the order of the convergence precision $10^{-12}$ for the growth scenario and below 10\% of the discretization error in the inflation scenario.

\vspace{2\baselineskip}

\subsection{Pinched Hemispherical Shell}

\begin{figure*}[htpb]
	\begin{center}
		\includegraphics{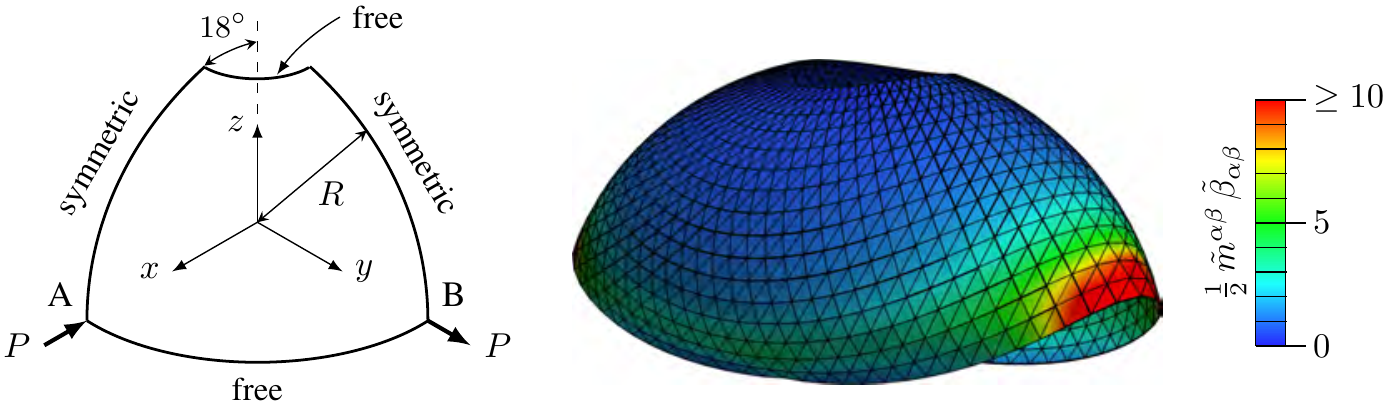}
		\caption{Hemispherical shell subject to point loads. Left: Undeformed reference configuration reduced to the first quadrant exploiting symmetry. Right: Control mesh of the quasi-static solution at maximum load $P=100$, without ghost elements. The color encodes the bending energy density.}
		\label{fig:hemisphere}
	\end{center}
\end{figure*}

Next, we turn to two numerical examples without growth, for verification of the linearly elastic, geometrically nonlinear shell with coupled stretching and bending. The pinched hemisphere is a widely used benchmark for ``an element's ability to represent inextensional modes'' and ``rigid body rotations about normals of the shell surface'' \cite{BLSCO85}. In its nonlinear regime, it is a test recommended by the National Agency for Finite Element Methods and Standards (NAFEMS, 3DNLG-9) \cite{ABM}. The geometrical setup is shown in Fig.~\ref{fig:hemisphere}. A hemispherical shell with an $18^\circ$ open pole and free boundaries is pinched by four equally strong, pairwise opposite diametrical point loads $P$ acting on the equator. The shell radius is $R=10$, its thickness $h=0.04$, and the elastic moduli are determined by $Y=6.825\!\times\!10^7$, $\nu=0.3$.

\begin{figure*}[htpb]
	\begin{center}
		\includegraphics{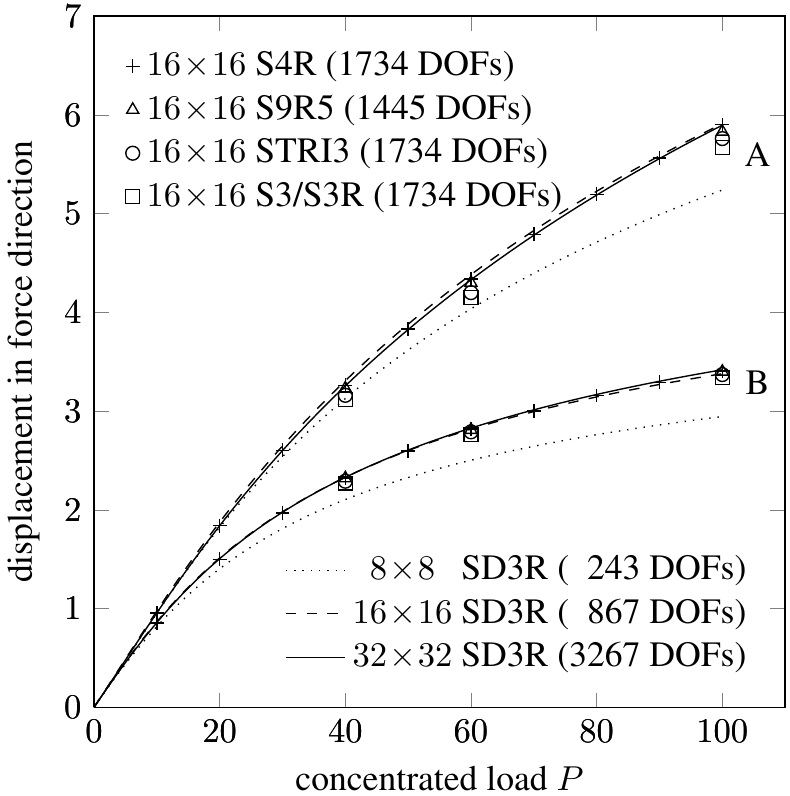}
		\caption{Load-displacement curves for the pinched hemispherical shell. The NAFEMS reference values were obtained with Abaqus.}
		\label{fig:hemisphere_plot}
	\end{center}
\end{figure*}

In order to minimize the impact originating from the specific choice of boundary constraints, the whole hemisphere is simulated without exploiting symmetry, and Green's method is used for the free boundaries, i.e., they are unconstrained. In Fig.~\ref{fig:hemisphere_plot}, we plot the displacements of points A and B on the limit surface against the applied loads for three mesh resolutions, together with some NAFEMS reference results obtained with Abaqus \cite{ABM,SLL04}. Much more reference data for other nonlinear finite elements can be found e.g.~in Refs.~\cite{ABM,NGMA09,W10} and references therein. The employed meshes are obtained by regularly discretizing the quarter hemisphere along the angles of inclination and the azimuth, resulting in 128, 512, and 2048 triangles per quarter hemisphere, respectively. In Fig.~\ref{fig:hemisphere}, the $16\!\times\!16$ mesh with 512 triangles per quarter, that is also used in the NAFEMS results, is shown on the right.

Our load-displacement curves for subdivision surface elements are almost identical to Abaqus S4R. High precision is obtained with SD3R at much less DOFs, mainly because subdivision shell elements go without rotational variables, unlike all other shown elements. The difference between one-point and six-point quadrature is well below 70\% of the discretization error in all cases.

\vspace{2\baselineskip}

\vspace{1\baselineskip}

\clearpage
\subsection{Stretched Cylinder with Free Ends}

\begin{figure*}[htpb]
	\begin{center}
		\includegraphics{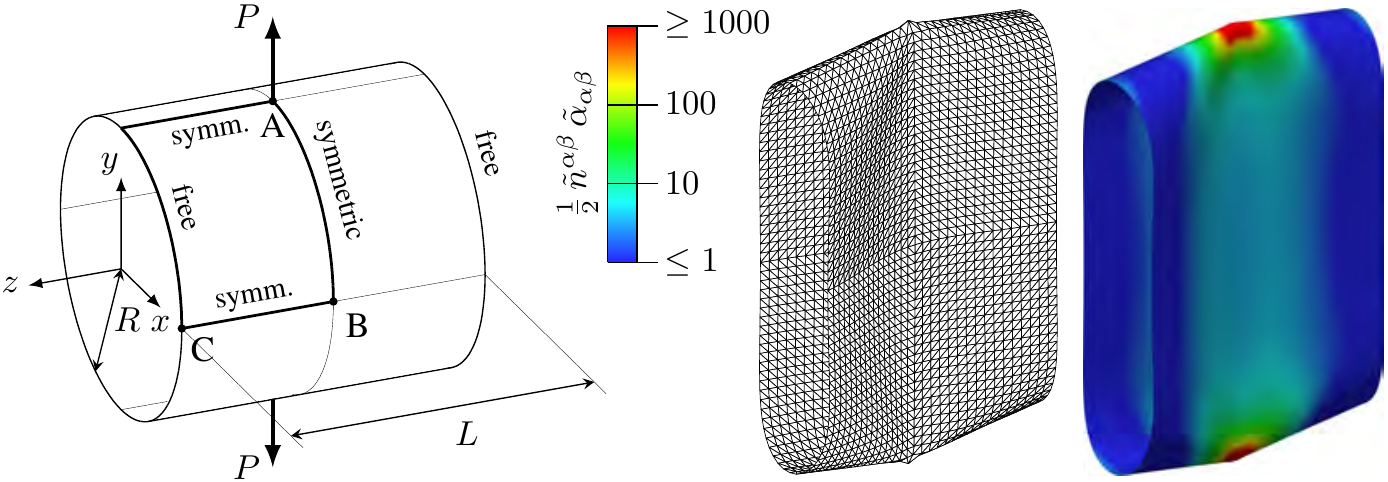}
		\caption{Cylindrical shell subject to point loads. Left: Undeformed reference configuration. Middle: Control mesh ($16\!\times\! 24$) of the quasi-static solution at maximum load $P=4\!\times\! 10^4$, without ghost elements. Right: Limit surface of the same configuration, with the stretching energy density on a logarithmic color scale.}
		\label{fig:cylinder}
	\end{center}
\end{figure*}

\begin{figure*}[htpb]
	\begin{center}
		\includegraphics{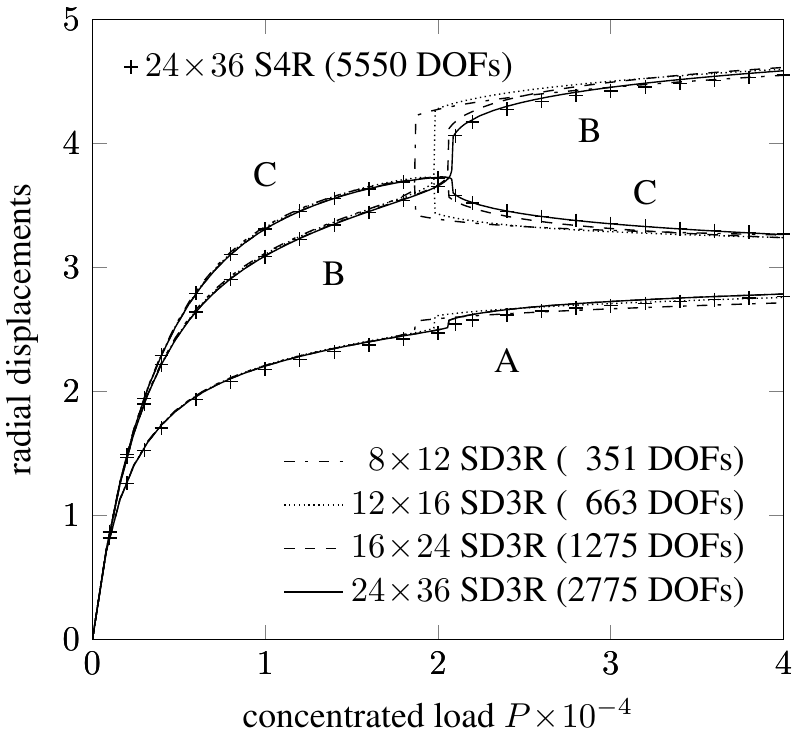}
		\caption{Load-displacement curves for the stretched cylindrical shell, with the transition near half loading.}
		\label{fig:cylinder_plot}
	\end{center}
\end{figure*}

The next example is another standard loading test, consisting of a cylindrical shell with free boundaries that is stretched transversally by two equally strong, opposite diametrical point loads $P$ acting on the middle of the cylinder length. This test case has found vast attention in the literature. For an overview, see eg.~Refs.~\cite{NGMA09,W10} and references therein. Its peculiar usefulness is due to its ability to examine two different response regimes, one after the other. At small loads, the large deformation results from low bending stiffness, while at large loads, further deformations require the stiff shell to be stretched primarily. The geometrical setup is shown in Fig.~\ref{fig:cylinder}. The cylinder radius is $R=4.953$, its length $L=10.35$, its thickness $h=0.094$, and the elastic moduli are determined by $Y=10.5\!\times\!10^6$, $\nu=0.3125$. Like in the previous example, the full shell is simulated neglecting present symmetries, and no constraints are applied to the boundaries. In Fig.~\ref{fig:cylinder_plot}, the resulting load-displacement curves for the points A, B and C are compared to Abaqus S4R element data from Ref.~\cite{SLL04}. In the bending regime at moderate loads, our data from 351 DOFs almost coincides with S4R using 5550 DOFs, again demonstrating the outstanding computational efficiency of subdivision shells. In the stretching regime at large loads, SD3R elements are slightly less stiff than S4R and elements found in other literature mentioned above. Those displacements, however, are quite sensitive to changes in the mesh structure. Other meshes than that shown in Fig.~\ref{fig:cylinder} lead to marginally shifted displacements at large loads. The relative error between one-point and six-point quadrature is at most equal to and mostly well below 40\% of the discretization error.

Correctly capturing the snap-through transition near $P\approx 2\!\times\!10^4$ has posed a tough challenge to various finite shell elements in the past. Some even fail to correctly feature it at moderate mesh resolutions \cite{GSW89,SB92,BDP95}. With the subdivision shell elements, we have not observed such problems, not even for meshes much coarser than mentioned in Fig.~\ref{fig:cylinder_plot}.

\vspace{3\baselineskip}

\FloatBarrier
\subsection{Anisotropic Growth and Boundary Instabilities}

\begin{figure*}[htpb]
	\begin{center}
		\includegraphics{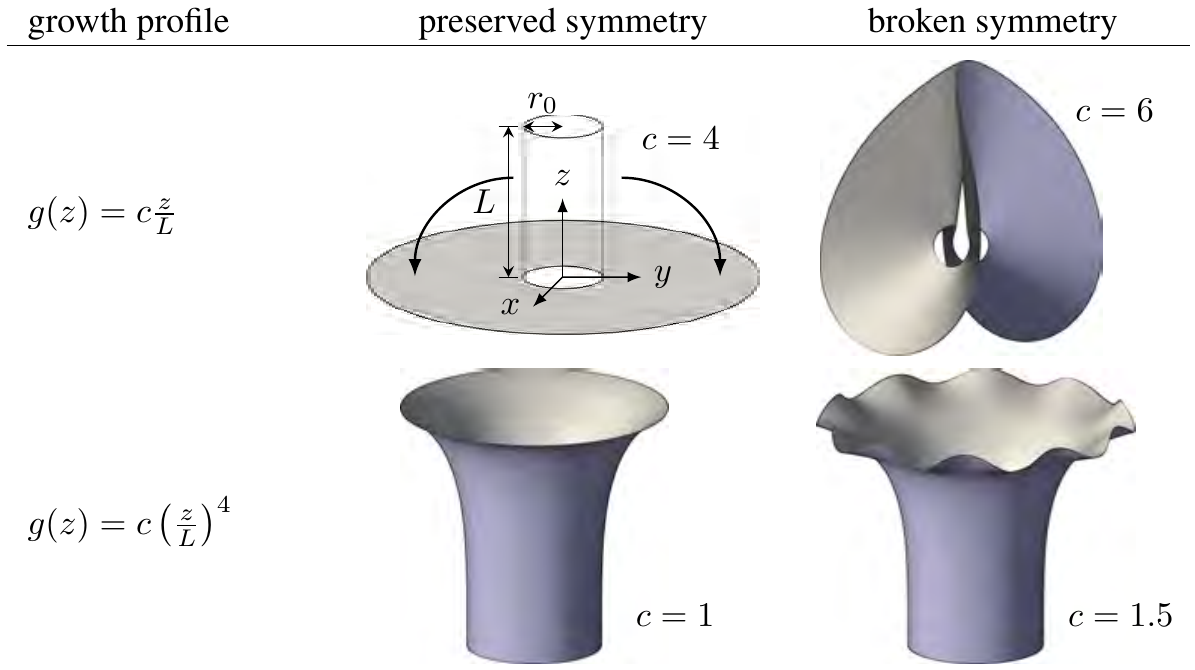}
		\caption{Symmetry breaking for polynomial growth gradients. The reference configuration is cylindrical with radius $r_0=1$ and length $L=4$. All boundaries are free, and the shell thickness is $h=10^{-2}$. The linear case with 1.5 times the critical growth leads to a ring (ground state: wavenumber $k=2$) with contact, while the same overcritical growth in the quartic case yields a $k=10$ mode at the boundary.}
		\label{fig:cylinder_growth}
	\end{center}
\end{figure*}

The purpose of the remaining two numerical examples is to further verify the proper functioning of the growing shell elements, and to demonstrate the wide range of possible applications. The first, presented in this subsection, addresses an interesting property of plant growth: the development of buckling instabilities at tissue boundaries, such as the edges of leaves and flower petals. The phenomenon is very similar to the permanent deformation left on a plastic sheet that is torn apart. Along the open boundary, wavy structures, whose exact nature depends on the shell thickness and growth profile, occur. The out-of-plane bending is a result of the tissue's thinness combined with large in-plane strains from compression due to in-plane growth (or plastic tension in the tearing case). When a certain critical growth threshold is reached, it becomes favorable to bend rather than compress further, and the initial symmetry is broken. Consider a shell whose reference configuration is a circular cylinder with radius $r(z)\equiv r_0$, $z\in[0,L]$, and assume an anisotropic growth tensor
\begin{equation}
{\bf G} = \textrm{diag}\big(1,1+g(z),1\big)
\end{equation}
with respect to the canonical basis of cylindrical coordinates $(r,\varphi,z)$. The growth profile $g(z)$ is considered monotonically increasing. Following the arguments in Ref.~\cite{M02}, where the Gauss\textendash{}Bonnet theorem is applied, the stability condition for cylindrical symmetry in the vanishing thickness limit reads
\begin{equation}
\left|r_0\frac{\textrm{d}g}{\textrm{d}z}(L)\right| \le 1.
\end{equation}
We have simulated polynomial growth fields of the form
\begin{equation}
g(z) = c\left(\frac{z}{L}\right)^{p},\qquad c,p > 0,
\end{equation}
for illustration and verification, where the cylindrical symmetry is preserved if $c\le L/r_0p$. A selection of fundamental results is shown in Fig.~\ref{fig:cylinder_growth}. The boundaries are of the free Green type, i.e., the ghost node displacements are unconstrained. The essential result of these simulations is that for reasonable angular resolutions of the meshed cylinder, the symmetry transition is observed very close to the theoretical value. Notice that the linear case $p=1$ is equivalent to the \textit{excess cone} \cite{DBA08, MBA08, SWAMH10} in the limit $r_0\rightarrow 0$, where a second instability is known to exist due to self-contact.

\begin{figure*}[h]
	\begin{center}
		\includegraphics{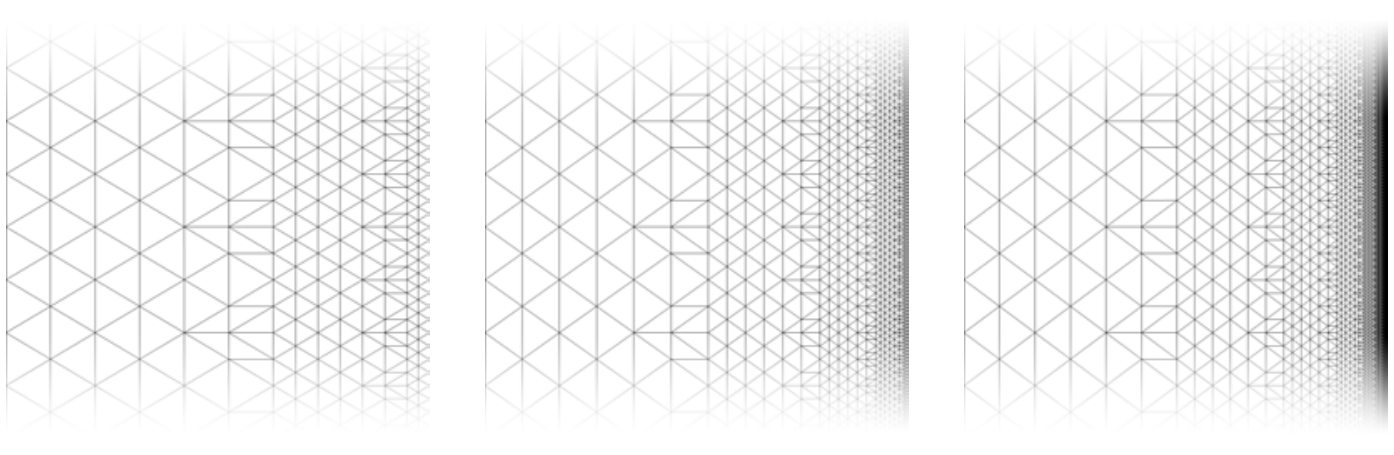}
		\caption{Employed mesh layouts. The ghost element layers are not shown.}
		\label{fig:mesh}
	\end{center}
\end{figure*}

\begin{figure*}[ht]
	\begin{center}
		\includegraphics{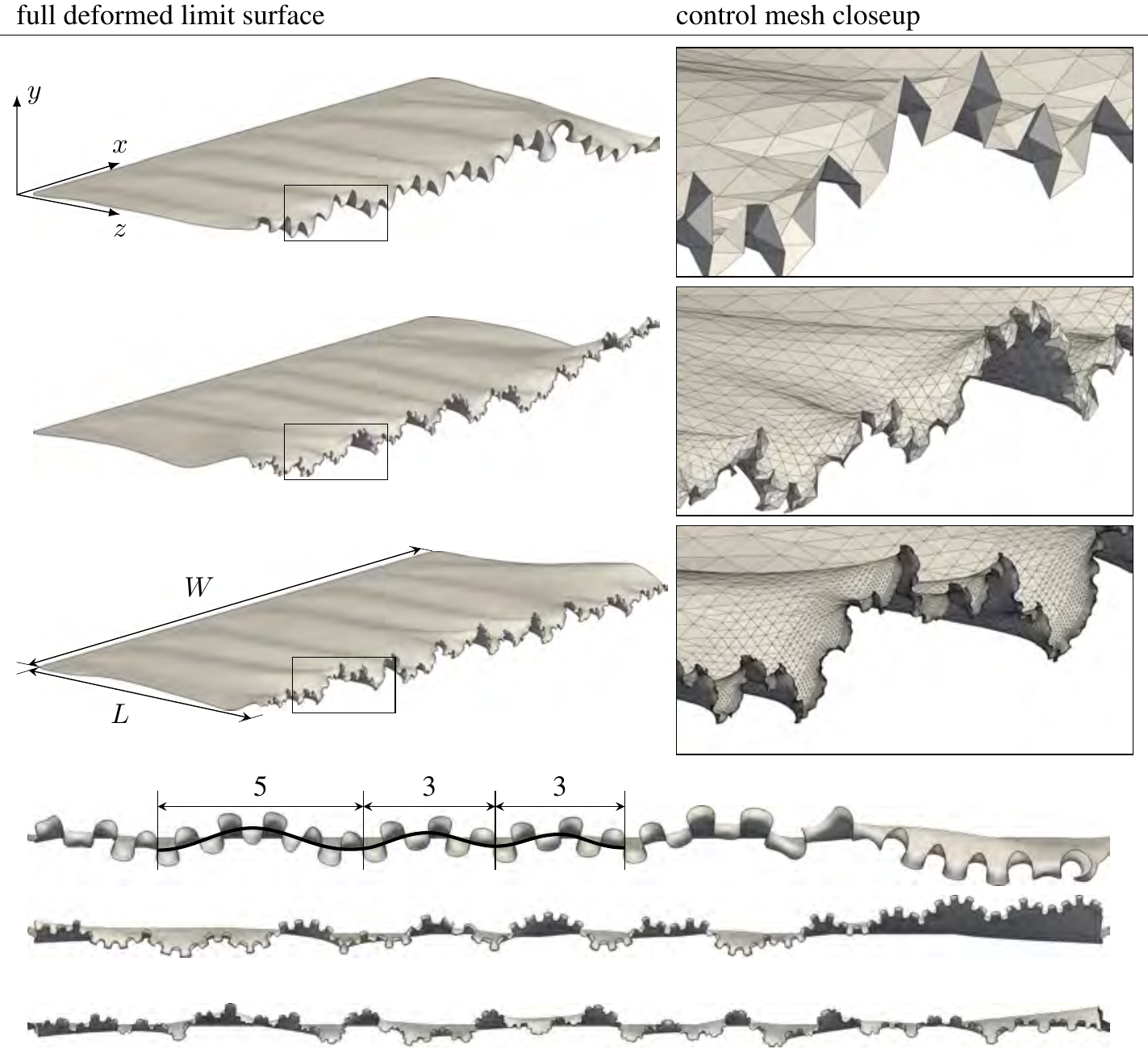}
		\caption{Self-similar shell boundary at large growth gradients on different mesh resolutions. All edges are free except at $z=0$, where the sheet is clamped in all directions. The ghost element layers are not shown. In the coarsest configuration (top) already, superimposed waves with scaling factor 3 and 5, as predicted in Ref.~\cite{AB03}, are observable. The finest resolution (bottom) is detailed enough to manifest about three levels of wrinkles, clearly resembling experimentally obtained self-similar wrinkling cascades \cite{SRMSS02,MDS07}.}
		\label{fig:selfsim_growth}
	\end{center}
\end{figure*}

For growth profiles with very large gradients $(\textrm{d}g/\textrm{d}z)(L)$ at the boundary, such as
\begin{equation}
g(z) = \left(1+\frac{L-z}{l}\right)^{-1},\qquad 0 < l \ll L,
\end{equation}
where $l$ is a small characteristic length scale, the pattern of the open shell boundary has been reported to be self-similar, with an odd integer scaling factor that is mostly 3 and sometimes close to 5 in experiments \cite{SRMSS02,AB03,A04}. This behavior is easily reproducible with our growth implementation. To demonstrate this, we have simulated a planar rectangular shell of length $L=1$, width $W=4$ and thickness $h=10^{-4}$, using a characteristic growth length $l=40h$. The $z=0$ edge is clamped using Schweitzer-Cirak constraints, while all other edges are unconstrained. The shell is grown using three different mesh resolutions to reveal the discretization effects. The employed hierarchical meshes are visualized in Fig.~\ref{fig:mesh}. Detailed mesh data is listed in Tab.~\ref{tab:mesh}, where also the $z$-coordinates of the Gauss points closest to the open edge, at which the growth tensors are largest, are given for reproducibility. Fig.~\ref{fig:selfsim_growth} shows the resulting equilibrated deformed configurations after growth. The self-similar nature of the open boundary is discernible even on the coarsest mesh, but the finest resolution is needed to get a clear resemblance to experiments \cite{SRMSS02,MDS07}.

\begin{table*}[!htbp]
	\centering
	\begin{tabular}{rrrrrr}
	\hline
	nodes & ghost nodes & elements & ghost elements & shortest element & largest quadrature $z$\\
	\hline
	1034 & 186 & 1876 & 376 & 1/32 & 0.999457\\
	4702 & 590 & 8804 & 1188 & 1/128 & 0.999887\\
	38780 & 4204 & 73340 & 8422 & 1/1024 & 0.999987\\
	\hline
	\end{tabular}
	\caption{Mesh properties for Figs.~\ref{fig:mesh} and \ref{fig:selfsim_growth}.}
	\label{tab:mesh}
\end{table*}

\setlength\fboxsep{0pt}

\subsection{Confined Growth vs.~Crumpling}

We put the growing shell to a final sophisticated test by comparing two processes involving tight self-interaction and spatial confinement, in the frictionless, quasi-static, elastic limit. A thin disk with radius $R$ is placed inside of a spherical container with the same radius ${\overline R}=R$. In the first setup, the container is shrunken, and consequently, the disk crumples into a ball of the size of the container, similarly to a sheet of paper that one crumples by hand. In the second setup, the opposite happens, i.e., the container sustains its size while the disk is subjected to a constant isotropic growth rate, both in plane and in thickness.

Various numerical simulations of crumpled elastic \cite{KW97,VG06,TAT08} and elasto-plastic \cite{TAT09} sheets and membranes in shrinking spheres have been carried out in recent years. The main finding is that sheets tightly crumpled into balls, although consisting mostly of air, develop a very large bulk stiffness resulting from a network of ridges and vertices of high magnitudes of mean curvature. A very large portion of the bending energy is condensed into this network \cite{LGLMW95}. A priori, it is not obvious whether a shell growing inside of a fixed confinement will exhibit equivalent behavior. With the present thin shell theory, we are able to answer this question in the elastic limit.

\begin{figure*}[htpb]
	\begin{center}
		\includegraphics{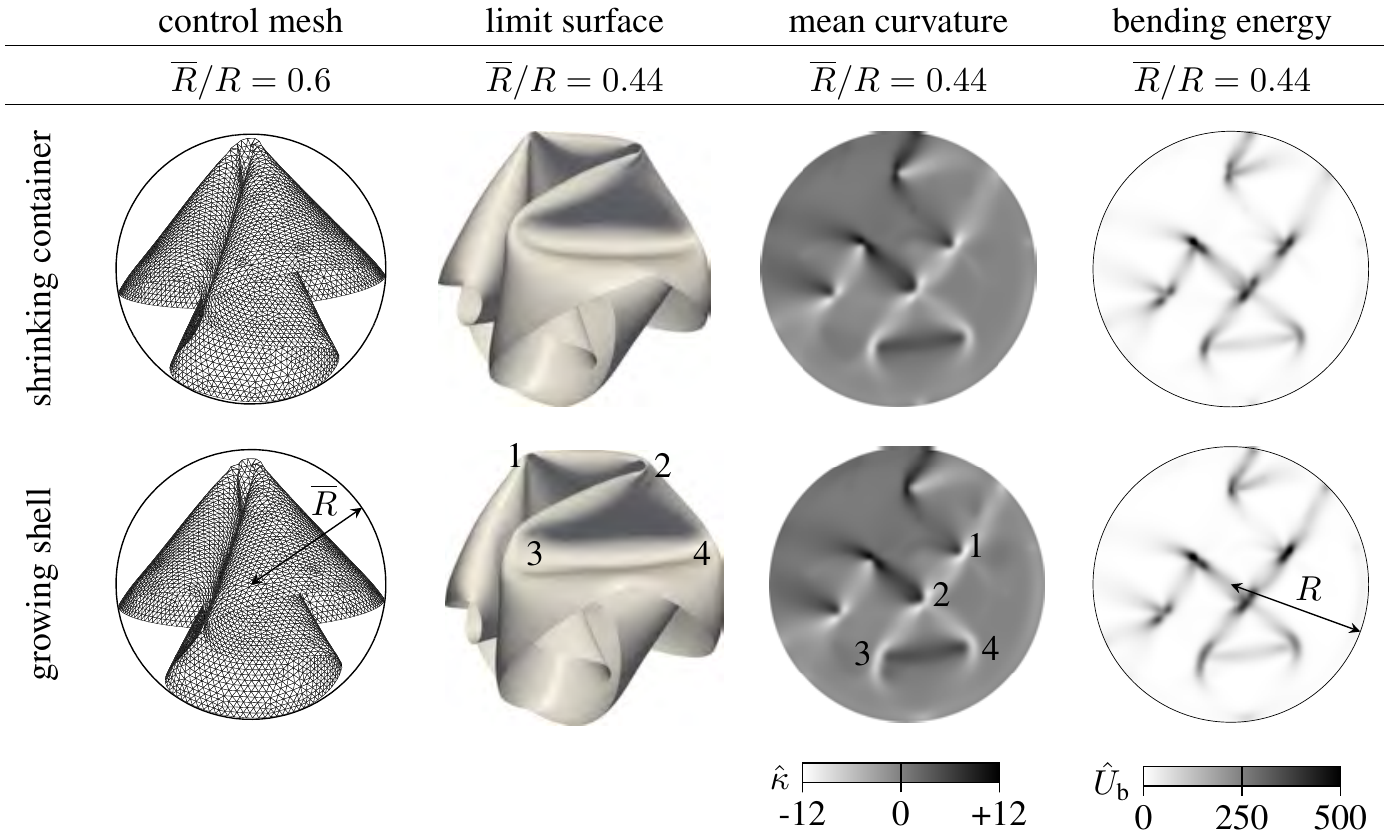}
		\caption{Comparison between a crumpled shell (top row) and a growing shell in spherical confinement (bottom row). Network vertices are numbered for easy identification.}
		\label{fig:shrink_grow}
	\end{center}
\end{figure*}

\clearpage
The only relevant simulation parameter is the thickness-to-size ratio $h/\overline{R}=0.01$, yielding a F\"oppl\textendash von K\'arm\'an number of
\begin{equation}\gamma=12(1-\nu^2)(\overline{R}/h)^2\approx 10^5.
\end{equation}
To obtain equivalent time scales in the two problems, we shrink the sphere in the first case according to $\overline{R}(t)=R/(1+g(t))$, where $g(t)=\lambda t$ is the in-plane growth factor on the growing shell in the second case. Accordingly, to achieve equivalent length scales, the growing shell has an increasing thickness $h(t)=h(0)(1+g(t))$. The growth rate $\lambda$ is chosen small enough to allow for a quasi-static simulation in both cases, and damping is subcritical. The used mesh consists of 8864 triangles (excluding ghosts, see Fig.~\ref{fig:shrink_grow}, first column).

We have not found any evidence indicating that the two processes are different. During early stages, both shells buckle to form a \textit{developable cone} with a single vertex, where most curvature and bending energy is concentrated. Around ${\overline R}/R\approx 0.53$, the apex splits into two vertices, and more vertices subsequently emerge, leading to the same ridge network. In the third column of Fig.~\ref{fig:shrink_grow}, the reduced mean curvature
\begin{equation}
\hat{\kappa} := \frac{k_1 + k_2}{2}{\overline R}
\end{equation}
is projected onto the unfolded disk, where $k_1$ and $k_2$ are the principal curvatures. The cross correlation of the mean curvature ridge patterns is $r=0.89$, a very high value when compared to recent similar measurements \cite{TAT08}. The fourth column displays the dimensionless rescaled bending energy density
\begin{equation}
\hat{U}_\textrm{b} := H^{\alpha\beta\gamma\delta}\tilde{\beta}_{\alpha\beta} \tilde{\beta}_{\gamma\delta}\overline{R}^2\left(\frac{R(t)}{R(0)}\right)^2.
\end{equation}
No qualitative or quantitative disparity, going beyond minor local shifts resulting from the finite element size and slow but finite dynamics, is observed. While this indicates that the two physical processes are in fact equivalent, it also accentuates the large potential of numerical simulations of growing thin shells.

\vspace{3\baselineskip}

\section{Conclusions}

We have presented a thin shell finite element approach that is capable of calculating anisotropic in-plane growth and large deformations. To this end, the Rodriguez deformation gradient decomposition was used to extend the classical Kirchhoff-Love theory. Assuming that the growth tensor is independent of the current deformation, the implementation is straightforward and requires only minor modifications, as it formally coincides with the implementation of Kirchhoff-Love shells without growth. The presented model generalizes the recently popularized target metric approach in that it includes the change of reference curvature of the shell.

The presented finite elements are of the Loop subdivision surface kind. Subdivision surface shape functions provide many advantages over traditional $C^0$ elements and other $C^1$ shape functions:
\begin{itemize}
\setlength{\itemsep}{0pt}
\setlength{\parskip}{0pt}
\item It meets the continuity requirement imposed by the bending energy and thus permits a shell finite element description in the classical Rayleigh\textendash{}Ritz formalism with all its amenities like optimal convergence.
\item One Gauss point per element is sufficient for Kirchhoff-Love shells.
\item They are easy to implement and require only three DOFs per node, nine per element. No rotational or auxiliary variables are needed.
\item They can handle arbitrary mesh topologies.
\item Subdivision surface elements are easily reformulated for other curvature-based theories. In particular, their applicability in the context of lipid bilayer mechanics was recently demonstrated \cite{FK06,MK08,KSM12a,KSM12b}.
\end{itemize}
We have demonstrated the outstanding efficiency of the subdivision shell elements on three standard loading examples. Despite their clear superiority to $C^0$ elements in this regard, subdivision shell elements are accompanied by new challenges:
\begin{itemize}
\setlength{\itemsep}{0pt}
\setlength{\parskip}{0pt}
\item The extended local support of subdivision basis functions requires new concepts for boundary conditions, adaptive mesh refinement and fracture modeling \cite{GKS02,G032,KTZ04,COP05}, which are more intricate than conventional techniques applicable to $C^0$ finite elements or discrete elements (e.g.~, \cite{WKHK04}).
\item Contact detection on the limit surface is a nonlinear optimization problem.
\end{itemize}

The finite element method is much better suited for strong material anisotropies than many other discretization schemes. In this article, anisotropic in-plane growth fields have been built into thin shell finite elements without complicating the underlying formalisms. We have illustrated the large range of potential applicability of growing finite shell elements to various problems in material science and engineering by numerically simulating different growth scenarios. They are also expected to be very well suited for the simulation of deformable confining membranes in packing problems \cite{SNWHH11,VWSH13}. A dry friction model and plasticity effects will have to be included to simulate real-world time-irreversible contact problems.

\vspace{\baselineskip}
\thanks{The authors acknowledge support from the ETH Research Grants ``Morphogenesis in Constrained Spaces'' and ``Packing of Slender Objects in Deformable Confinements'' (ETHIIRA Grants No. TH-06 07-3 and ETH-03 10-3) and from the European Research Council (ERC) Advanced Grant 319968-FlowCCS. V.~Lienhard has assisted in implementing shell-shell contact handling.}


\begin{thebibliography}{99}
\setlength{\itemsep}{-1ex}

\bibitem{O98}
Otubushin A. Detailed validation of a non-linear finite element code using
  dynamic axial crushing of a square tube. \emph{International Journal of
  Impact Engineering}  1998; \textbf{21}(5):349--368.

\bibitem{WWK01}
Webb DC, Webster J, Kormi K. Finite {E}lement {S}imulation of {E}nergy
  {A}bsorption {D}evices under {A}xial {S}tatic {C}ompressive and {I}mpact
  {L}oading. \emph{International Journal of Crashworthiness}  2001;
  \textbf{6}(3):399--424.

\bibitem{MMIKH01}
Mamalis A, Manolakos D, Ioannidis M, Kostazos P, Hassiotis G. Finite element
  simulation of the axial collapse of thin-wall square frusta.
  \emph{International Journal of Crashworthiness}  2001;
  \textbf{6}(2):155--164.

\bibitem{MO05}
Maia LG, de~Oliveira PHIA. A {R}eview of {F}inite {E}lement {S}imulation of
  {A}ircraft {C}rashworthiness. \emph{SAE Technical Paper}  2005.

\bibitem{LGLMW95}
Lobkovsky A, Gentges S, Li H, Morse D, Witten TA. Scaling properties of
  stretching ridges in a crumpled elastic sheet. \emph{Science}  1995;
  \textbf{270}(5241):1482--1485.

\bibitem{BAP97}
Ben~Amar M, Pomeau Y. Crumpled paper. \emph{Proceedings of the Royal Society A}
   1997; \textbf{453}(1959):729--755.

\bibitem{BK05}
Blair DL, Kudrolli A. Geometry of {C}rumpled {P}aper. \emph{Physical Review
  Letters}  2005; \textbf{94}:166107.

\bibitem{VG06}
Vliegenthart GA, Gompper G. Forced crumpling of self-avoiding elastic sheets.
  \emph{Nature Materials}  2006; \textbf{5}:216--221.

\bibitem{TAT08}
Tallinen T, \AA{}str\"om JA, Timonen J. Deterministic {F}olding in {S}tiff
  {E}lastic {M}embranes. \emph{Physical Review Letters}  2008;
  \textbf{101}:106101.

\bibitem{TAT09}
Tallinen T, \AA{}str\"om JA, Timonen J. The effect of plasticity in crumpling
  of thin sheets. \emph{Nature Materials}  2009; \textbf{8}:25--28.

\bibitem{VG11}
Vliegenthart GA, Gompper G. Compression, crumpling and collapse of spherical
  shells and capsules. \emph{New Journal of Physics}  2011;
  \textbf{13}(4):045020.

\bibitem{L88}
Love AEH. The {S}mall {F}ree {V}ibrations and {D}eformation of a {T}hin
  {E}lastic {S}hell. \emph{Philosophical Transactions of the Royal Society of
  London. (A.)}  1888; \textbf{179}:491--546.

\bibitem{P95}
Parisch H. A continuum-based shell theory for non-linear applications.
  \emph{International Journal for Numerical Methods in Engineering}  1995;
  \textbf{38}(11):1855--1883.

\bibitem{ZT05}
Zienkiewicz OC, Taylor RL. \emph{The Finite Element Method for Solid and
  Structural Dynamics}. 6th edn., Elsevier Butterworth-Heinemann: Oxford, 2005.
  ISBN 0-7506-6321-9.

\bibitem{SC07}
Stogner RH, Carey GF. ${C}^1$ macroelements in adaptive finite element methods.
  \emph{International Journal for Numerical Methods in Engineering}  2007;
  \textbf{70}(9):1076--1095.

\bibitem{COS00}
Cirak F, Ortiz M, Schr\"oder P. Subdivision surfaces: a new paradigm for
  thin-shell finite-element analysis. \emph{International Journal for Numerical
  Methods in Engineering}  2000; \textbf{47}(12):2039--2072.

\bibitem{CO01}
Cirak F, Ortiz M. Fully ${C}^1$-conforming subdivision elements for finite
  deformation thin-shell analysis. \emph{International Journal for Numerical
  Methods in Engineering}  2001; \textbf{51}(7):813--833.

\bibitem{T95}
Taber LA. Biomechanics of {G}rowth, {R}emodeling, and {M}orphogenesis.
  \emph{Applied Mechanics Reviews}  1995; \textbf{48}(8):487--545.

\bibitem{DBA08}
Dervaux J, Ben~Amar M. Morphogenesis of {G}rowing {S}oft {T}issues.
  \emph{Physical Review Letters}  2008; \textbf{101}:068101.

\bibitem{H86}
Hoger A. On the determination of residual stress in an elastic body.
  \emph{Journal of Elasticity}  1986; \textbf{16}:303--324.

\bibitem{SZJNH96}
Skalak R, Zargaryan S, Jain RK, Netti PA, Hoger A. Compatibility and the
  genesis of residual stress by volumetric growth. \emph{Journal of
  Mathematical Biology}  1996; \textbf{34}:889--914.

\bibitem{KHBSH12}
Kim J, Hanna JA, Byun M, Santangelo CD, Hayward RC. Designing responsive
  buckled surfaces by halftone gel lithography. \emph{Science}  2012;
  \textbf{335}(6073):1201--1205.

\bibitem{KHHS12}
Kim J, Hanna JA, Hayward RC, Santangelo CD. Thermally responsive rolling of
  thin gel strips with discrete variations in swelling. \emph{Soft Matter}
  2012; \textbf{8}:2375--2381.

\bibitem{H68}
Hsu FH. The influences of mechanical loads on the form of a growing elastic
  body. \emph{Journal of Biomechanics}  1968; \textbf{1}(4):303--311.

\bibitem{RHM94}
Rodriguez EK, Hoger A, McCulloch AD. Stress-dependent finite growth in soft
  elastic tissues. \emph{Journal of Biomechanics}  1994;
  \textbf{27}(4):455--467.

\bibitem{DQ02}
DiCarlo A, Quiligotti S. Growth and balance. \emph{Mechanics Research
  Communications}  2002; \textbf{29}(6):449--456.

\bibitem{LH02}
Lubarda V, Hoger A. On the mechanics of solids with a growing mass.
  \emph{International Journal of Solids and Structures}  2002;
  \textbf{39}(18):4627--4664.

\bibitem{L69}
Lee EH. Elastic-{P}lastic {D}eformation at {F}inite {S}trains. \emph{ASME
  Journal of Applied Mechanics}  1969; \textbf{36}(1):1--6.
  
\bibitem{AG05}
Ambrosi D, Guana F. Stress-modulated growth. \emph{Mathematics and Mechanics of
  Solids}  2007; \textbf{12}(3):319--342.

\bibitem{MG11}
Moulton D, Goriely A. Anticavitation and {D}ifferential {G}rowth in {E}lastic
  {S}hells. \emph{Journal of Elasticity}  2011; \textbf{102}:117--132.

\bibitem{K66}
Koiter WT. On the nonlinear theory of thin elastic shells. \emph{Proceedings of
  the Koninklijke Nederlandse Akademie van Wetenschappen Royal Dutch Academy of
  Sciences Ser B}  1966; \textbf{69}:1--54.

\bibitem{C00}
Ciarlet P. Un mod\`ele bi-dimensionnel non lin\'eaire de coque analogue \`a
  celui de {W}.{T}.~{K}oiter. \emph{Comptes Rendus de l'Acad\'emie des Sciences
  - Series I - Mathematics}  2000; \textbf{331}(5):405--410.

\bibitem{ESK09}
Efrati E, Sharon E, Kupferman R. Elastic theory of unconstrained
  non-{E}uclidean plates. \emph{Journal of the Mechanics and Physics of Solids}
   2009; \textbf{57}(4):762--775.

\bibitem{SF89}
Simo J, Fox D. On a stress resultant geometrically exact shell model. {P}art
  {I}: {F}ormulation and optimal parametrization. \emph{Computer Methods in
  Applied Mechanics and Engineering}  1989; \textbf{72}(3):267--304.

\bibitem{SFR90}
Simo J, Fox D, Rifai M. On a stress resultant geometrically exact shell model.
  {P}art {III}: {C}omputational aspects of the nonlinear theory. \emph{Computer
  Methods in Applied Mechanics and Engineering}  1990; \textbf{79}(1):21--70.

\bibitem{G03}
Green S. Multilevel, subdivision-based, thin shell finite elements: Development
  and an application to red blood cell modelings. Ph{D} {T}hesis, University of
  Washington, Seattle 2003.

\bibitem{MSSR03}
Marder M, Sharon E, Smith S, Roman B. Theory of edges of leaves.
\emph{Europhysics Letters}  2003; \textbf{62}(4):498.

\bibitem{AB03}
Audoly B, Boudaoud A. Self-similar structures near boundaries in strained
  systems. \emph{Physical Review Letters}  2003; \textbf{91}:086105.

\bibitem{KES07}
Klein Y, Efrati E, Sharon E. Shaping of elastic sheets by prescription of
  non-euclidean metrics. \emph{Science}  2007; \textbf{315}(5815):1116--1120.

\bibitem{MDS07}
Marder M, Deega RD, Sharon E. {C}rumpling, {B}uckling, and {C}rackling:
  {E}lasticity of {T}hin {S}heets. \emph{Physics Today}  2007;
  \textbf{60}:33--38.

\bibitem{CC78}
Catmull E, Clark J. Recursively generated b-spline surfaces on arbitrary
  topological meshes. \emph{Computer-Aided Design}  1978;
  \textbf{10}(6):350--355.

\bibitem{S99}
Stam J. Evaluation of {L}oop {S}ubdivision {S}urfaces. \emph{SIGGRAPH '99
  course notes}, 1999.

\bibitem{L87}
Loop C. Smooth subdivision surfaces based on triangles. Master's {T}hesis,
  Department of Mathematics, University of Utah, Salt Lake City 1987.

\bibitem{S71}
Stroud AH. \emph{{A}pproximate {C}alculation of {M}ultiple {I}ntegrals}.
  Prentice\textendash{}Hall: Englewood Cliffs, 1971. ISBN 0-13-043893-6.

\bibitem{CL11}
Cirak F, Long Q. Subdivision shells with exact boundary control and
  non-manifold geometry. \emph{International Journal for Numerical Methods in
  Engineering}  2011; \textbf{88}(9):897--923.

\bibitem{BLZ00}
Biermann H, Levin A, Zorin D. Piecewise smooth subdivision surfaces with normal
  control. \emph{Proceedings of the 27th annual conference on Computer graphics
  and interactive techniques}, SIGGRAPH '00, 2000; 113--120.

\bibitem{S96}
Schweitzer JE. Analysis and {A}pplication of {S}ubdivision {S}urfaces. Ph{D}
  {T}hesis, Department of Computer Science and Engineering, University of
  Washington, Seattle 1996.

\bibitem{G04}
Green S, Turkiyyah G. Second-order accurate constraint formulation for
  subdivision finite element simulation of thin shells. \emph{International
  Journal for Numerical Methods in Engineering}  2004; \textbf{61}(3):380--405.

\bibitem{SFR89}
Simo J, Fox D, Rifai M. On a stress resultant geometrically exact shell model.
  {P}art {II}: {T}he linear theory; {C}omputational aspects. \emph{Computer
  Methods in Applied Mechanics and Engineering}  1989; \textbf{73}(1):53--92.

\bibitem{N59}
Newmark NM. A {M}ethod of {C}omputation for {S}tructural {D}ynamics.
  \emph{Journal of the Engineering Mechanics Division}  1959; \textbf{85}(EM
  3):67--94.

\bibitem{ZWLX92}
Zeng LF, Wiberg NE, Li XD, Xie YM. \textit{A posteriori} local error estimation
  and adaptive time-stepping for {N}ewmark integration in dynamic analysis.
  \emph{Earthquake Engineering \& Structural Dynamics}  1992;
  \textbf{21}(7):555--571.

\bibitem{R04}
Reddy JN. \emph{An Introduction to Nonlinear Finite Element Analysis}. Oxford
  University Press: Oxford, 2004; 292--295. ISBN 0-19-852529-X.

\bibitem{ZX91}
Zienkiewicz OC, Xie YM. A simple error estimator and adaptive time stepping
  procedure for dynamic analysis. \emph{Earthquake Engineering \& Structural
  Dynamics}  1991; \textbf{20}(9):871--887.

\bibitem{GCSO99}
Grinspun E, Cirak F, Schr\"oder P, Ortiz M. {N}on-{L}inear {M}echanics and
  {C}ollisions for {S}ubdivision {S}urfaces. \emph{Technical {R}eport},
  California Institute of Technology, Pasadena 1999.

\bibitem{R60}
Rosen JB. The {G}radient {P}rojection {M}ethod for {N}onlinear {P}rogramming.
  {P}art {I}. {L}inear {C}onstraints. \emph{Journal of the Society for
  Industrial and Applied Mathematics}  1960; \textbf{8}(1):181--217.

\bibitem{BFA02}
Bridson R, Fedkiw R, Anderson J. Robust treatment of collisions, contact and
  friction for cloth animation. \emph{ACM Transactions on Graphics}  2002;
  \textbf{21}(3):594--603.

\bibitem{E05}
Eriscon C. \emph{Real-Time Collision Detection}. Morgan Kaufmann: San
  Francisco, 2005; 136--142. ISBN 1-55860-732-3.

\bibitem{S01}
Sunday D. Distance between {L}ines and {S}egments with their {C}losest {P}oint
  of {A}pproach 2001.
  \\\url{http://www.softsurfer.com/Archive/algorithm_0106/algorithm_0106.htm}.

\bibitem{AUM}
Dassault Syst\`emes Simulia Corp., Providence, RI, USA. \emph{Abaqus 6.11
  {U}ser's {M}anual} 2011.

\bibitem{CB98}
Chapelle D, Bathe KJ. Fundamental considerations for the finite element
  analysis of shell structures. \emph{Computers \& Structures}  1998;
  \textbf{66}(1):19--36.

\bibitem{CB11}
Chapelle D, Bathe KJ. \emph{The Finite Element Analysis of Shells \textendash{}
  Fundamentals}. 2nd edn., Springer, 2011. ISBN 978-3-642-16407-1.

\bibitem{BIC00}
Bathe KJ, Iosilevich A, Chapelle D. An evaluation of the {MITC} shell elements.
  \emph{Computers \& Structures}  2000; \textbf{75}(1):1--30.

\bibitem{GS50}
Green AE, Shield RT. Finite elastic deformation of incompressible isotropic
  bodies. \emph{Proceedings of the Royal Society of London. Series A.
  Mathematical and Physical Sciences}  1950; \textbf{202}(1070):407--419.

\bibitem{N77}
Needleman A. Inflation of spherical rubber balloons. \emph{International
  Journal of Solids and Structures}  1977; \textbf{13}(5):409--421.

\bibitem{BLSCO85}
Belytschko T, Stolarski H, Liu WK, Carpenter N, Ong JS. Stress projection for
  membrane and shear locking in shell finite elements. \emph{Computer Methods
  in Applied Mechanics and Engineering}  1985; \textbf{51}(1-3):221--258.

\bibitem{ABM}
Dassault Syst\`emes Simulia Corp., Providence, RI, USA. \emph{Abaqus 6.11
  {B}enchmarks {M}anual} 2011.

\bibitem{SLL04}
Sze K, Liu X, Lo S. Popular benchmark problems for geometric nonlinear analysis
  of shells. \emph{Finite Elements in Analysis and Design}  2004;
  \textbf{40}(11):1551--1569.

\bibitem{NGMA09}
Negahban M, Goel A, Marchon P, Azizinamini A. Geometrically {E}xact {N}onlinear
  {E}xtended-{R}eissner/{M}indlin {S}hells: {F}undamentals, {F}inite {E}lement
  {F}ormulation, {E}lasticity. \emph{International Journal for Computational
  Methods in Engineering Science and Mechanics}  2009; \textbf{10}(6):430--449.

\bibitem{W10}
Wi\`{s}niewski K. \emph{Finite Rotation Shells: Basic Equations and Finite
  Elements for Reissner Kinematics}. Springer Netherlands: Dordrecht, 2010.
  ISBN 90-481-8760-5.

\bibitem{GSW89}
Gruttmann F, Stein E, Wriggers P. Theory and numerics of thin elastic shells
  with finite rotations. \emph{Archive of Applied Mechanics}  1989;
  \textbf{59}:54--67.

\bibitem{SB92}
Sansour C, Bufler H. An exact finite rotation shell theory, its mixed
  variational formulation and its finite element implementation.
  \emph{International Journal for Numerical Methods in Engineering}  1992;
  \textbf{34}(1):73--115.

\bibitem{BDP95}
Brank B, Damjani\'{c} FB, Peri\'{c} D. On implementation of a nonlinear four
  node shell finite element for thin multilayered elastic shells.
  \emph{Computational Mechanics}  1995; \textbf{16}:341--359.

\bibitem{M02}
Marder M. The {S}hape of the {E}dge of a {L}eaf  2002; ArXiv:cond-mat/0208232.

\bibitem{MBA08}
M\"uller MM, Ben~Amar M, Guven J. Conical {D}efects in {G}rowing {S}heets.
  \emph{Physical Review Letters}  2008; \textbf{101}:156104.

\bibitem{SWAMH10}
Stoop N, Wittel FK, Ben~Amar M, M\"uller MM, Herrmann HJ. Self-{C}ontact and
  {I}nstabilities in the {A}nisotropic {G}rowth of {E}lastic {M}embranes.
  \emph{Physical Review Letters}  2010; \textbf{105}:068101.

\bibitem{SRMSS02}
Sharon E, Roman B, Marder M, Shin GS, Swinney HL. Mechanics: {B}uckling
  cascades in free sheets. \emph{Nature (London)}  2002; \textbf{419}:579.

\bibitem{A04}
Audoly B. The self-similar rippling of leaf edges and torn plastic sheets.
  \emph{Europhysics News}  2004; \textbf{35}(5):145--148.

\bibitem{KW97}
Kramer EM, Witten TA. Stress condensation in crushed elastic manifolds.
  \emph{Physical Review Letters}  1997; \textbf{78}:1303--1306.

\bibitem{FK06}
Feng F, Klug WS. Finite element modeling of lipid bilayer membranes.
  \emph{Journal of Computational Physics}  2006; \textbf{220}(1):394--408.

\bibitem{MK08}
Ma L, Klug WS. Viscous regularization and r-adaptive remeshing for finite
  element analysis of lipid membrane mechanics. \emph{Journal of Computational
  Physics}  2008; \textbf{227}(11):5816--5835.

\bibitem{KSM12a}
Kahraman O, Stoop N, M\"uller MM. Morphogenesis of membrane invaginations in
  spherical confinement. \emph{Europhysics Letters}  2012;
  \textbf{97}(6):68008.

\bibitem{KSM12b}
Kahraman O, Stoop N, M\"uller MM. Fluid membrane vesicles in confinement.
  \emph{New Journal of Physics}  2012; \textbf{14}(9):095021.

\bibitem{GKS02}
Grinspun E, Krysl P, Schr\"oder P. {CHARMS}: {A} {S}imple {F}ramework
  for {A}daptive {S}imulation. \emph{SIGGRAPH (ACM Transactions on Graphics)}
  2002; \textbf{21}(3):281--290.

\bibitem{G032}
Grinspun E. {T}he {B}asis {R}efinement {M}ethod. Ph{D} {T}hesis, Department of
  Computer Science, Columbia University, New York 2003.

\bibitem{KTZ04}
Krysl P, Trivedi A, Zhu B. Object-oriented hierarchical mesh refinement with
  {CHARMS}. \emph{International Journal for Numerical Methods in Engineering}
  2004; \textbf{60}(8):1401--1424.

\bibitem{COP05}
Cirak F, Ortiz M, Pandolfi A. A cohesive approach to thin-shell fracture and
  fragmentation. \emph{Computer Methods in Applied Mechanics and Engineering}
  2005; \textbf{194}(21-24):2604--2618.

\bibitem{WKHK04}
Wittel F, Kun F, Herrmann HJ, Kr\"oplin BH. Fragmentation of {S}hells.
  \emph{Physical Review Letters}  2004; \textbf{93}:035504.

\bibitem{SNWHH11}
Stoop N, Najafi J, Wittel FK, Habibi M, Herrmann HJ. {P}acking of {E}lastic
  {W}ires in {S}pherical {C}avities. \emph{Physical Review Letters}  2011;
  \textbf{106}:214102.

\bibitem{VWSH13}
Vetter R, Wittel FK, Stoop N, Herrmann HJ. Finite element simulation of dense
  wire packings. \emph{European Journal of Mechanics - A/Solids}  2013;
  \textbf{37}:160--171.

\end{thebibliography}
\end{document}